\documentclass[a4paper,11pt,leqno,notitlepage,twoside]{amsart}

\usepackage[T1]{fontenc}
\usepackage[latin1]{inputenc}
\usepackage[frenchb]{babel}
\usepackage{amsmath,amssymb,amsthm,amsfonts,url}
\usepackage{geometry}
\usepackage{epic}
\usepackage{graphicx,color}
\usepackage[right]{eurosym}
\usepackage{ifpdf}
\usepackage{array}
\usepackage{marvosym}
\usepackage{wasysym}

\numberwithin{equation}{section}
  
\theoremstyle{plain} 
\newtheorem{theorem}{Théorème}[section] 
\newtheorem{corollary}[theorem]{Corollaire}

\theoremstyle{definition} 
\newtheorem{definition}{Définition}[section]
\newtheorem{remark}{Remarque}[section]

\title[Hauteur asymptotique des points de Heegner]{Hauteur asymptotique des points de Heegner} 
\author[G. Ricotta - T. Vidick]{Guillaume Ricotta - Thomas Vidick}
\thanks{2005 \textit{Mathematics Subject Classification:} 11G50, 11M41.}
\date{Version of \today}

\newcommand{\mn}{\mathnormal}
\newcommand{\mf}{\mathfrak}
\newcommand{\Z}{\mathbb{Z}}
\newcommand{\N}{\mathbb{N}}

\newcommand{\C}{\mathbb{C}}
\newcommand{\Q}{\mathbb{Q}}

\newcommand{\f}[2]{{\ensuremath{\mathchoice%
        {\dfrac{#1}{#2}}
        {\dfrac{#1}{#2}}
        {\frac{#1}{#2}}
        {\frac{#1}{#2}}
        }}}
\newcommand{\paf}[2]{\ensuremath{\left(\f{#1}{#2}\right)}}
\newcommand{\pa}[1]{\ensuremath{\left(#1\right)}}

\begin{document}


\renewcommand{\abstractname}{Abstract}
\begin{abstract}
The asymptotic behaviour of the N\'eron-Tate height of Heegner points on a rational elliptic curve attached to an arithmetically normalized new cusp form $f$ of weight $2$, level $N$ and trivial character is studied in this paper. By Gross-Zagier formula, this height is related to the special value at the critical point for the derivative of the Rankin-Selberg convolution of $f$ with a certain weight one theta series attached to some ideal class of some imaginary quadratic field. Asymptotic formula for the first moments asociated to these Dirichlet series are proved and experimental results are carefully discussed.  
\end{abstract}

\maketitle

\noindent{\textit{En l'honneur du Professeur Henryk Iwaniec, pour l'ensemble de son oeuvre analytique et pour son titre de Docteur Honoris Causa de l'Université de Bordeaux 1.}}

\tableofcontents

\section{Description de la problématique}

\begin{figure}[htbp]
\begin{center}
\includegraphics[height=10cm]{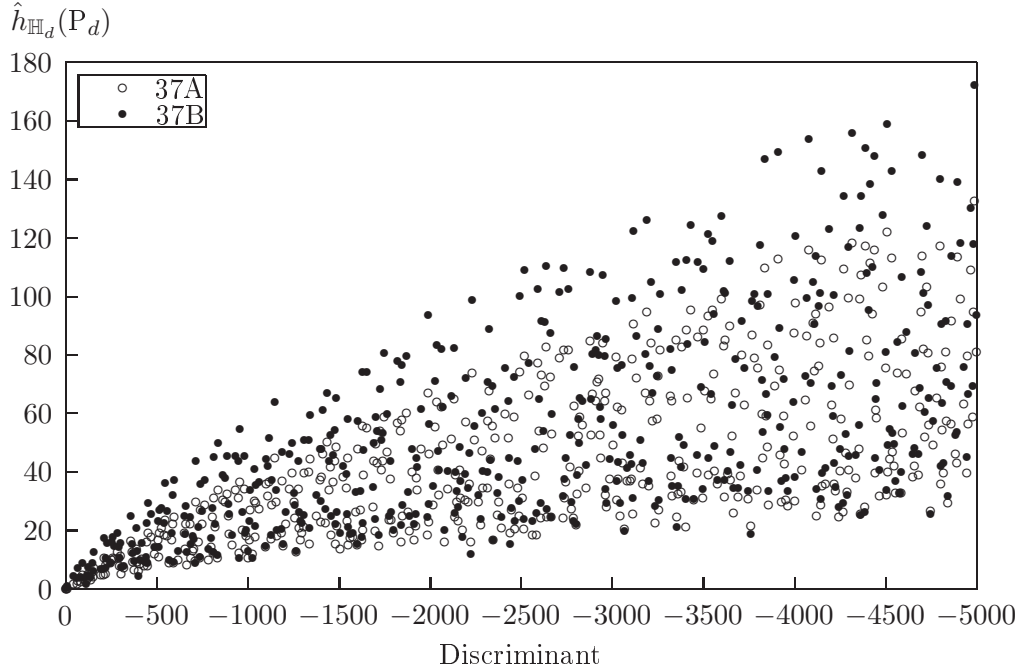}
\end{center}
\caption{Comportement des points de Heegner sur les courbes $37\text{A}$ et $37\text{B}$}
\label{37A-1}
\end{figure}

L'\'etude calculatoire syst\'ematique de la hauteur de N\'eron-Tate des points de Heegner sur diff\'erentes courbes elliptiques montre de grandes disparit\'es. Si l'on consid\`ere deux courbes elliptiques de m\^eme conducteur $N$ alors les points de Heegner sur ces deux courbes sont l'image des m\^emes points spéciaux de la courbe modulaire de niveau $N$ notée $\text{X}_0(N)$ par la param\'etrisation modulaire et on s'attend donc \`a ce que leurs hauteurs soient dans le m\^eme rapport que les degr\'es de ces param\'etrisations. Cependant, la figure \ref{37A-1} montre que le comportement des hauteurs est tr\`es irr\'egulier: même si les courbes $37\text{A}$ et $37\text{B}$ (dans la notation de Cremona (\cite{cremona})) ont le m\^eme conducteur et le m\^eme degr\'e, les points paraissent l\'eg\`erement plus gros sur la $37\text{B}$ que sur la $37\text{A}$. 

Nous allons montrer que l'intuition se v\'erifie asymptotiquement: \`a conducteur fix\'e\footnote{Toutefois, nous prendrons soin de garder explicite toute dépendance en le conducteur $N$ de la courbe.}, la moyenne sur une certaine sous-classe de discriminants de la hauteur des points de Heegner est proportionnelle au degr\'e.

Pour cela, nous allons proc\'eder \`a partir de la formule de Gross-Zagier (\cite{GZ}) et raisonner sur les s\'eries de Dirichlet et les fonctions $L$ en nous servant principalement d'un r\'esultat de H. Iwaniec (\cite{iwaniec}). 

\noindent{\textbf{Remerciements}.} Les auteurs remercient chaleureusement H. Darmon pour leur avoir suggéré cette problématique et pour ses nombreux conseils et encouragements. Ce travail a été réalisé à l'occasion d'un stage d'études à McGill University (Montréal) pour le second auteur et d'un stage post-doctoral à l'Université de Montréal (Montréal) pour le premier auteur. Les excellentes conditions de travail offertes par ces deux institutions ont fortement contribué à la réalisation de cet article. Le premier auteur a largement profité de la générosité et des conseils Mathématiques avisés de A. Granville lors de son stage post-doctoral.

\section{Pr\'eliminaires}

Soit $E$ une courbe elliptique d\'efinie sur $\Q$. Supposons pour simplifier que son conducteur $N$ est sans facteurs carr\'es et considérons l'ensemble
$$\mathcal{D}:= \mn{\{d\in\mathbb{Z}_-^*, \mu^2(d)=1, d\equiv \nu^2 \mod 4N, (\nu,4N)=1\}}$$
de discriminants. Pour $d$ dans $\mathcal{D}$, notons $\mathbb{H}_d$ le corps de classe de Hilbert du corps quadratique imaginaire $\mathbb{K}_d:=\Q(\sqrt{d})$ et souvenons-nous que $G_d:=\text{Gal}(\mathbb{H}_d\vert\mathbb{K}_d)$ est isomorphe au groupe de classes de $\mathbb{K}_d$ de cardinal le nombre de classes $h_d$. Soit $X_0(N)$ la courbe modulaire de niveau $N$  classifiant les paires de courbes elliptiques $(E_1,E_2)$ reli\'ees par une isog\'enie cyclique de degr\'e $N$. Une description analytique sur $\mathbb{C}$ de cette courbe est donnée par le quotient du demi-plan de Poincaré complété $\mathbb{H}\cup\left(\mathbb{Q}\cup\left\{\infty\right\}\right)$ par l'action par homographies du groupe de congruence $\varGamma_0(N)$. 
\begin{definition}
Un \emph{point de Heegner de niveau $N$ et de discriminant $d$} est un couple ordonné $(E_1,E_2)$ de courbes elliptiques muni d'une isogénie cyclique de degré $N$ tel que $E_1$ et $E_2$ aient multiplication complexe par l'anneau des entiers $\mathcal{O}_d$ de $\mathbb{K}_d$. 
\end{definition}
Fixons une fois pour toute une racine carr\'ee $s_d$ de $d$ modulo $4N$ et désignons par $\mf{n}_d$ l'id\'eal entier primitif de norme $N$ suivant:
$$\mf{n}_d\mn{:=\pa{N,\f{s_d+\sqrt{d}}{2}}}.$$
À $s_d$ fixé, l'ensemble des points de Heegner de niveau $N$ et de discriminant $d$ est en bijection avec le groupe de classes de $\mathbb{K}_d$ au sens suivant: si $[\mathfrak{a}]$ est l'élément du groupe de classes de $\mathbb{K}_d$ associé à l'idéal entier primitif $\mathfrak{a}$ de $\mathcal{O}_d$ alors $(\C/\mf{a}\mn{,\C/}\mf{a}\mathfrak{n}_d{^{-1}})$ est un point de Heegner de niveau $N$ et de discriminant $d$. Le point du demi-plan de Poincaré modulo $\varGamma_0(N)$ correspondant à ce point de Heegner est donné par $\frac{-B+\sqrt{d}}{2A}$ modulo $\varGamma_0(N)$ où $(A,B,C)$ est la forme quadratique de discriminant $d$ correspondant à $[\mathfrak{a}]$ et où $N\mid A$ et $B\equiv s_d\mod 2N$.
\begin{definition}
Un \emph{point de Heegner de niveau $N$ et de discriminant $d$ sur $E$} est l'image par la param\'etrisation modulaire $\Phi_{N,E}: X_0(N)\rightarrow E$ d'un point de la forme
$$(E_1,E_2)=(\C/\mf{a}\mn{,\C/}\mf{a}\mathfrak{n}_d{^{-1}}),$$
o\`u $\mf{a}$ est un id\'eal de l'anneau des entiers $\mathcal{O}_d$ de $\mathbb{K}_d$.
\end{definition}
Notons P$_d=\Phi_{N,E}\left(\left(E_1,E_2\right)\right)$, ce point ne d\'epend que de la classe de $\mf{a}$ dans le groupe de classes de $\mathbb{K}_d$ (une fois $E$ et $d$ fix\'es). Notons \'egalement Tr$_d=$Tr$_{\mathbb{H}_d\vert\mathbb{K}_d}($P$_d)$. Les points de Heegner sont d\'efinis sur $E(\mathbb{H}_d)$ et sont permut\'es par $G_d$ (\cite{gross}). \newline\newline
Si $R$ est un corps de nombres alors $\hat{h}_R(P)$ désigne la hauteur de N\'eron-Tate (comme d\'efinie dans \cite{si86}, VIII.9) prise sur $R$ du point $P$ à coordonnées dans $R$. On rappelle que si $S$ est une extension de degr\'e fini de $R$ alors $\hat{h}_S(P)=[S:R]\hat{h}_R(P)$. Le but de cet article est d'\'etudier la valeur en moyenne, sur les $d$ dans $\mathcal{D}$ des deux objets suivants:
\begin{itemize}
\item 
$\hat{h}_{\mathbb{H}_d}($P$_d)$, hauteur de N\'eron-Tate sur $\mathbb{H}_d$ d'un quelconque des $h_d$ points de Heegner d\'efinis ci-dessus. Puisque $\hat{h}_{\mathbb{H}_d}$ est invariante sous l'action de $G_d$, cette hauteur est ind\'ependante du point choisi,
\item 
$\hat{h}_{\mathbb{K}_d}($Tr$_d)$, hauteur de N\'eron-Tate de la trace sur $\mathbb{K}_d$ d'un quelconque des points de Heegner P$_d$ d\'efinis ci-dessus.
\end{itemize}
B.H. Gross et D. Zagier (\cite{GZ}) ont reli\'e $\hat{h}_{\mathbb{H}_d}($P$_d)$ \`a la valeur de la d\'eriv\'ee en $1$ de la série de Dirichlet obtenue en effectuant le produit de la fonction $L$ de Dirichlet associée au caractère primitif réel $\chi_d(m):=\paf{d}{m}$ de conducteur $\vert d\vert$ du corps $\mathbb{K}_d$ (le caractère de Kronecker du corps) par la convolution de Rankin-Selberg de $L(E\vert\mathbb{Q},s)$ avec la fonction zeta $\sum_{n\geqslant 1}r_d(n)n^{-s}$ de la classe des idéaux principaux de $\mathbb{K}_d$ c'est-à-dire
\begin{align}
\label{lesd}
L_d(E,s):=\left(\sum_{\substack{m\geqslant 1 \\
(m,N)=1}}\frac{\chi_d(m)}{m^{2s-1}}\right)\times\left(\sum_{n\geqslant 1}\f{a_n r_d(n)}{n^s}\right)
\end{align}
o\`u pour tout entier naturel non-nul $n$, $r_d(n)$ d\'esigne le nombre d'id\'eaux principaux de $\mathbb{K}_d$ de norme $n$.
\begin{theorem}[B.H. Gross-D. Zagier (1986)]
Si $E$ est une courbe elliptique d\'efinie sur $\Q$ et $d$ est dans $\mathcal{D}$ alors
\begin{align}
\label{gzpoint}
L^\prime_d(E,1)= \f{2\Omega_{E,N}}{u^2\sqrt{-d}} \hat{h}_{\mathbb{H}_d}(\emph{P}_d),\\
\label{gztrace}
L^\prime(E\vert\mathbb{K}_d,1)=\f{2\Omega_{E,N}}{u^2\sqrt{-d}} \hat{h}_{\mathbb{K}_d}(\emph{Tr}_d)
\end{align}
o\`u $L(E\vert\mathbb{K}_d,s)$ est la fonction $L$ de $E$ sur le corps $\mathbb{K}_d$, $2u$ est le nombre de racines de l'unit\'e de $\mathbb{K}_d$ et $\Omega_{E,N}=Im(\omega_1\bar{\omega}_2)$ est le volume complexe de $E$ (c'est-à-dire le double de l'aire d'un parallélogramme fondamental de $E(\mathbb{C})$).
\end{theorem}
Nous rappelons finalement la relation \'evidente 
$$\hat{h}_{\mathbb{K}_d}(\text{Tr}_d)=\hat{h}_{\mathbb{H}_d}(\text{P}_d)+\sum_{\sigma\in G_d\backslash \{Id\}} <\text{P}_d,\text{P}_d^{\sigma}>_{\mathbb{H}_d}$$
entre nos deux objets d'\'etude ($<\cdot,\cdot>_{\mathbb{H}_d}$ d\'esigne la forme bilin\'eaire de N\'eron-Tate sur $\mathbb{H}_d$). Le troisi\`eme terme repr\'esente l'angle form\'e par les points de Heegner entre eux et il sera int\'eressant de l'analyser (voir paragraphe \ref{5}).

\section{Les traces}
\label{traces}

\subsection{Mise en place}

Il est ici n\'ecessaire de raisonner selon le rang de la courbe $E$. En effet, soit $L(E\vert\mathbb{Q},s):=\sum_{n\geqslant 1}a_n n^{-s}$ sa série $L$ sur $\Q$ définie a priori sur $\Re{(s)}>\frac{3}{2}$. Selon les travaux de A. Wiles et de R. Taylor (\textbf{\cite{Wi}} et \textbf{\cite{TaWi}}), il existe une forme primitive cuspidale $f$ de niveau $N$, de poids $2$ et de caractère trivial telle que
\begin{equation*}
L(E\vert\mathbb{Q},s)=L(f,s)
\end{equation*}
sur $\Re{(s)}>\frac{3}{2}$. Par conséquent, $L(E\vert\mathbb{Q},s)$ admet un prolongement holomorphe à $\mathbb{C}$ et satisfait l'équation fonctionnelle
\begin{equation*}
\left(\frac{\sqrt{N}}{2\pi}\right)^s\varGamma(s)L(E\vert\mathbb{Q},s)=\omega\left(\frac{\sqrt{N}}{2\pi}\right)^{2-s}\varGamma(2-s)L(E\vert\mathbb{Q},2-s)
\end{equation*}
où $\omega=\pm 1$ est une valeur propre d'Atkin-Lehner de $f$. Définissons pour tout discriminant $d$ dans $\mathcal{D}$ la fontion $L$ de $E$ sur $\mathbb{Q}$ tordue par $\chi_d$ sur $\Re{(s)}>\frac{3}{2}$ par
$$L(E\vert\mathbb{Q}\times\chi_d,s):=\sum_{n\geqslant 1}\frac{a_n\chi_d(n)}{n^s}.$$
$L(E\vert\mathbb{Q}\times\chi_d,s)$ admet un prolongement holomorphe à $\mathbb{C}$ et satisfait l'équation fonctionnelle
\begin{equation*}
\left(\frac{\vert d\vert\sqrt{N}}{2\pi}\right)^s\varGamma(s)L(E\vert\mathbb{Q}\times\chi_d,s)=\omega_d\left(\frac{\vert d\vert\sqrt{N}}{2\pi}\right)^{2-s}\varGamma(2-s)L(E\vert\mathbb{Q}\times\chi_d,2-s)
\end{equation*}
où $\omega_d:=\omega\chi_d(-N)=-\omega$ (voir \cite{IwKo}). Avec ces notations, la factorisation
\begin{equation}
\label{lesurk}
L(E\vert\mathbb{K}_d,s)=L(E\vert\mathbb{Q},s)L(E\vert\mathbb{Q}\times\chi_d,s)
\end{equation}
est valide (voir \cite{Da}) d'où
\begin{equation*}
L^\prime(E\vert\mathbb{K}_d,1)=L^\prime(E\vert\mathbb{Q},1)L(E\vert\mathbb{Q}\times\chi_d,1)+L(E\vert\mathbb{Q},1)L^\prime(E\vert\mathbb{Q}\times\chi_d,1).
\end{equation*}
Ainsi, l'\'etude de la hauteur des traces des points de Heegner en moyenne sur les discriminants $d$ dans $\mathcal{D}$ est ramen\'ee \`a l'\'etude de la valeur au point critique des fonctions $L$ tordues lorsque $E$ est de rang analytique $1$ et \`a l'\'etude des d\'eriv\'ees au point critique des fonctions $L$ tordues lorsque $E$ est de rang analytique $0$ car alors $\omega=+1$. 

\subsection{Courbes de rang analytique $0$}

Rappelons le th\'eor\`eme obtenu par H. Iwaniec dans \cite{iwaniec}. Avant cela, fixons une fois pour toutes les deux notations
\begin{equation}
\label{gammaN}
\gamma(4N):=\#\{d\mod 4N, d\equiv \nu^2\mod 4N, (\nu,4N)=1\},
\end{equation}
et
\begin{equation}
\label{cN}
c_N:=\f{3\gamma(4N)}{\pi^2N}\prod_{\substack{p\in\mathcal{P} \\
p\mid 2N}} \pa{1-\f{1}{p^2}}^{-1}.
\end{equation}
\begin{theorem}
\label{iwa}
Si $E$ une courbe elliptique rationnelle de conducteur $N$ sans facteurs carrés et de rang analytique $\,\,0$ et $F$ est une fonction lisse à support compact dans $\mathbb{R}_+$ et de moyenne strictement positive alors 
$$\sum_{d\in\mathcal{D}} L'(E\vert\mathbb{Q}\times\chi_d,1)F\left(\frac{\vert d\vert}{Y}\right)=\alpha_NY\log{Y}+\beta_NY+\mathcal{O}_\varepsilon\left(N^{\frac{23}{14}+\varepsilon}Y^{\frac{13}{14}+\varepsilon}\right)$$ 
pour tout $\varepsilon>0$ o\`u
\begin{equation}
\label{alphaN}
\alpha_N:=c_NL(1)\int_0^{+\infty}F(t)\mathrm{d}t\neq 0
\end{equation}
et
\begin{equation*}
\beta_N:=c_N\int_0^{+\infty}F(t)\left(L^\prime(1)+L(1)\left(\log{\left(\frac{\sqrt{N}t}{2\pi}\right)}-\gamma\right)\right)\mathrm{d}t
\end{equation*}
avec
\begin{equation*}
L(s):=\f{L\left(\text{Sym}^2E,2s\right)}{\zeta^{(N)}(4s-2)}\left(\prod_{\substack{p\in\mathcal{P} \\
p\mid N}} \pa{1-\f{a_p}{p^s}}^{-1}\pa{1-\f{a_{p^2}}{p^{2s}}}\right)\mathcal{P}(s)
\end{equation*}
et
\begin{equation}
\label{P}
\mathcal{P}(s):=\prod_{\substack{p\in\mathcal{P} \\
p\nmid 2N}} \pa{\f{1}{1+1/p}+\pa{\f{1}{1+p}}\pa{\f{1+p^{2-4s}-(a_p^2-2p)p^{-2s}}{1+p^{1-2s}}}}.
\end{equation}
\end{theorem}
\begin{remark}
Dans \cite{iwaniec}, la fonction $L$ est définie par la série de Dirichlet suivante
$$L(s)=\sum_{\substack{n=k\ell^2 \\
k\mid N^\infty \\
(\ell,N)=1}}\frac{b_n}{n^{s}}$$
avec
$$b_n:=a_n\prod_{\substack{p\in\mathcal{P} \\
p\mid n \\
p\nmid 2N}} \pa{1+\f{1}{p}}^{-1}$$
pour tout entier naturel non-nul $n$. Celle-ci peut se réécrire sous la forme
\begin{equation*}
L(s)=\prod_{\substack{p\in\mathcal{P} \\
p\mid N}} \pa{1-\f{a_p}{p^s}}^{-1} \prod_{\substack{p\in\mathcal{P} \\
p\nmid 2N}} \pa{1+\pa{1+\f{1}{p}}^{-1}\pa{\sum_{i=1}^{\infty} \f{a_{p^{2i}}}{p^{2is}}}} \prod_{\substack{p\in\mathcal{P} \\
p\mid(2,N-1)}} \pa{1+\sum_{i=1}^{\infty} \f{a_{p^{2i}}}{p^{2is}}}
\end{equation*}
ce qui permet de retrouver l'expression de $L$ donnée dans le théorème en fonction du carr\'e sym\'etrique de $E$ et du produit Eulérien $\mathcal{P}$. Signalons que l'holomorphie de la fonction $L\left(\text{Sym}^2E,s\right)$ dans tout le plan complexe a \'et\'e prouv\'ee par Shimura dans \cite{shimura} et que le produit Eulérien $\mathcal{P}(s)$ est absolument convergeant sur $\Re{(s)}>\frac{3}{4}$ et y définit une fonction holomorphe.
\end{remark}
\begin{remark}
La valeur de la fonction $L\left(\text{Sym}^2E,s\right)$ au bord de la bande critique est reli\'ee au degr\'e de la param\'etrisation modulaire $\Phi_{N,E}:X_0(N)\rightarrow E$ de $E$ (cf. \cite{watkins} (1-1)) par la formule suivante analogue à celle du nombre de classes de Dirichlet
$$\f{L\left(\text{Sym}^2E,2\right)}{\pi\Omega_{E,N}}=\f{\text{deg}(\Phi_{N,E})}{Nc_E(N)^2}$$
où $c_E(N)$ est la $\varGamma_0(N)$-constante de Manin de $E$ qui est un entier relatif uniformément borné (\cite{Ed}). Ici, on utilise le fait que $N$ est sans facteurs carr\'es à deux reprises:
\begin{itemize}
\item
les fonctions $L$ du carré symétrique de $E$ motivique et analytique coincident car il n'y a pas de termes correctifs en les nombres premiers dont le carré divise le conducteur $N$,
\item
lorsque $E$ est une courbe de Weil $X_0(N)$-forte et $N$ est impair, la constante de Manin vaut $\pm 1$ selon les travaux de A. Abbes et E. Ullmo (\cite{AbUl}) alors que J. Manin (\cite{Ma}) a conjecturé que cette constante vaut $\pm 1$ pour toute courbe de Weil $X_0(N)$-forte\footnote{Ceci est faux si $E$ n'est pas forte: $[0,1,1,0,0]$ a une $X_0(11)$-constante de Manin égale à $5$.}.
\end{itemize}
\end{remark}
\begin{remark}
Il semblerait que quelques petites erreurs de frappe dans $\alpha_N$ (et en fait dans $c_N$ et $L$) se soient glissées dans \cite{iwaniec}. La valeur donnée ici est corrigée.
\end{remark}
\begin{remark}
Dans \cite{iwaniec}, la dépendance en le conducteur $N$ de la courbe dans le terme d'erreur n'est pas explicite. Cependant, il suffit de reprendre les différentes majorations pour restituer celle-ci.
\end{remark}
\begin{corollary}
\label{rang0}
Si $E$ est une courbe elliptique rationnelle de conducteur $N$ sans facteurs carrés et de rang analytique $0$ alors
\begin{multline}
\label{eqrg0}
\sum_{\substack{d\in\mathcal{D} \\
\vert d\vert\leqslant Y}}\hat{h}_{\mathbb{K}_d}(\emph{Tr}_d)=C_{\emph{Tr}}^{(0)}Y^{\frac{3}{2}}\log{Y}+\frac{C_{\emph{Tr}}^{(0)}}{2}\log{(N)}Y^{\frac{3}{2}} \\
+\frac{L(E\vert\mathbb{Q},1)c_N}{3\Omega_{E,N}}\left(L^\prime(1)-L(1)\left(\frac{2}{3}+\log{(2\pi)}+\gamma\right)\right)Y^{\frac{3}{2}}+\mathcal{O}_\varepsilon\left(\frac{L(E\vert\mathbb{Q},1)}{\Omega_{E,N}}N^{\frac{23}{14}+\varepsilon}Y^{\frac{20}{14}+\varepsilon}\right)
\end{multline}
pour tout $\varepsilon>0$ o\`u $C_{\emph{Tr}}^{(0)}$ est la constante d\'efinie par
$$C_{\emph{Tr}}^{(0)}:=\f{2}{\pi}c_N\mathcal{P}(1)L(E\vert\mathbb{Q},1)\f{L(\text{Sym}^2 E,2)}{\pi\Omega_{E,N}} L_{E}$$
en notant $L_{E}$ un produit de facteurs locaux correspondant aux facteurs Eulériens de $L(E\vert\mathbb{Q},s)$ et de $L(\text{Sym}^2(E),s)$ en les nombres premiers divisant le conducteur
\begin{equation}
\label{LE}
L_{E}:=\prod_{\substack{p\in\mathcal{P} \\
p\mid N}} \pa{1-\f{a_p}{p}}^{-1}\pa{1-\f{1}{p^2}}^{-1}\pa{1-\f{a_{p^2}}{p^2}}.
\end{equation}
\end{corollary}
\noindent{\textbf{Preuve du corollaire \ref{rang0}.}} Il suffit d'appliquer la formule de Gross-Zagier \eqref{gztrace} et de prendre pour $F$ une approximation lisse de la fonction qui vaut $\f{3}{2}\sqrt{t}$ sur $[0,1]$ et $0$ en dehors de cet intervalle. Le corollaire d\'ecoule alors de \eqref{lesurk}.
\begin{flushright}
$\blacksquare$
\end{flushright}
\begin{remark}
\`A conducteur $N$ fixé, le terme principal dans \eqref{eqrg0} est $C_{\emph{Tr}}^{(0)}Y^{\frac{3}{2}}\log Y$ avec
\begin{equation*}
C_{\emph{Tr}}^{(0)}=\left(\frac{6}{\pi^3c_E(N)^2}\mathcal{P}(1)\prod_{\substack{p\in \mathcal{P} \\
p \mid 2N}}\left(1-\frac{1}{p^2}\right)^{-1}\right)\times L(E\vert\mathbb{Q},1)\,L_{E}\times\frac{\gamma(4N)}{N^2}\text{deg}(\Phi_{N,E})
\end{equation*}
et est donc proportionnel au degré de la paramétrisation modulaire. Par contre, si $N=Y^a$ avec $0<a<\frac{1}{23}$ alors le deuxième terme de \eqref{eqrg0} est du même ordre de grandeur que le premier et le terme principal devient
\begin{equation*}
 \left(1+\frac{a}{2}\right)C_{\emph{Tr}}^{(0)}Y^{\frac{3}{2}}\log{Y}.
 \end{equation*}
C'est la principale raison pour laquelle nous avons rendu explicite la dépendance en le conducteur $N$ de la courbe dans le terme d'erreur. Il serait aussi intéressant d'étudier l'influence des valeurs extrémales de $L(E\vert\mathbb{Q},1)$ sur la moyenne des hauteurs des traces des points de Heegner. Malheureusement, cela ne semble pas vérifiable numériquement étant donné le temps de calcul nécessaire par les algorithmes à utiliser.
\end{remark}
Le produit Eulérien $\mathcal{P}(s)$ varie peu. En effet, si l'on effectue un développement limit\'e du facteur local en le nombre premier $p$ de $\mathcal{P}(1)$, on obtient
$$\mathcal{P}_{p}(1)=1-\f{a_p^2-2p}{p^3}+\mathcal{O}\pa{\f{1}{p^2}}$$
d'où l'existence d'une constante absolue $C>0$ ne dépendant pas de la courbe $E$ considérée telle que
$$C\prod_{p\in\mathcal{P}}\pa{1-\f{2}{p^2}}<\mathcal{P}(1)<C\prod_{p\in\mathcal{P}}\pa{1+\f{2}{p^2}.}$$
Finalement, les termes importants sont le degré de la paramétrisation modulaire, la valeur de la fonction $L$ de la courbe au point critique $1$ ainsi que le conducteur qui intervient en $\gamma(4N)/N^2$.
\begin{figure}[htbp]
\begin{center}
\includegraphics[height=10cm]{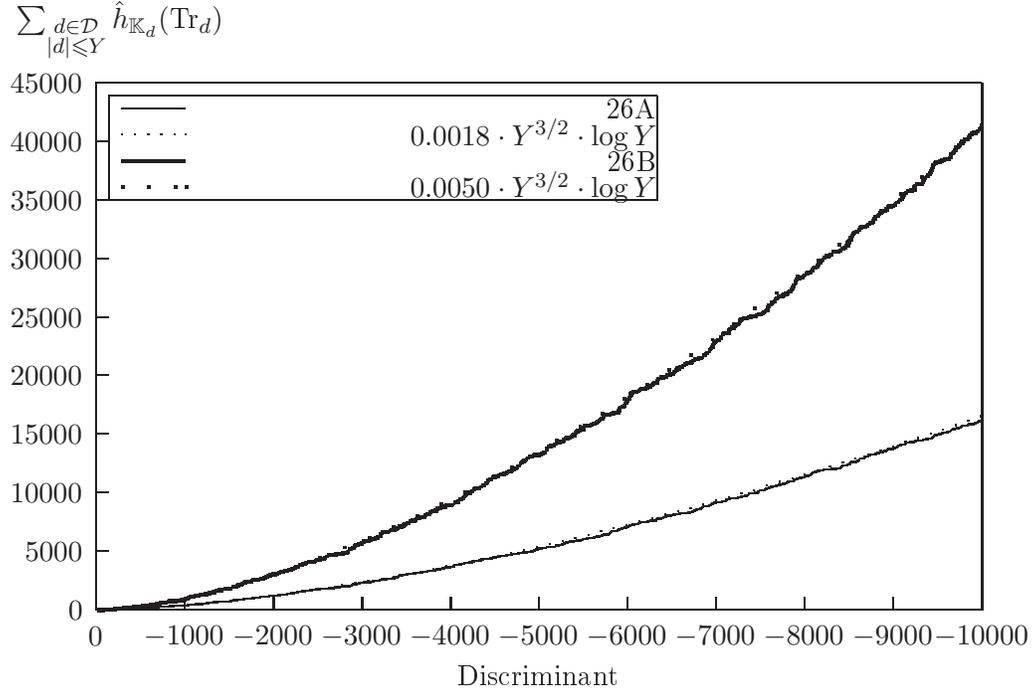}
\end{center}
\caption{Hauteur des traces en moyenne sur les courbes $26\text{A}$ et $26\text{B}$ : pratiques et th\'eoriques.}
\label{26A-traces}
\end{figure}
\`A conducteur et degré fixés, les petits nombres premiers qui divisent le conducteur ont un r\^ole décisif et les hauteurs des traces ont alors tendance à être plus ou moins grandes suivant que $E$ a réduction \emph{multiplicative d\'eploy\'ee} ou \emph{non-d\'eploy\'ee} en ces nombres premiers. Cette influence se voit très bien dans le cas des deux courbes de conducteur $26$ (de rang analytique $0$) dont les param\'etrisations modulaires ont même degré $2$ alors que les traces sont presque trois fois plus grosses sur la $26\text{B}$ que sur la $26\text{A}$. Ceci s'explique analytiquement par le fait que $26$ est divisible par $2$ et que la $26\text{A}$ a réduction multiplicative d\'eploy\'ee en $2$ alors que la $26\text{B}$ a r\'eduction multiplicative non-d\'eploy\'ee. Ceci induit un facteur $3$ entre les produits pour les nombres premiers $p$ divisant le conducteur de ces courbes des facteurs Euleriens de leur fonction $L$ en $p$. La figure \ref{26A-traces} illustre cette diff\'erence de comportement. On voit \'egalement que la courbe repr\'esentant la somme des hauteurs des traces est tr\`es proche de la courbe th\'eorique m\^eme si elle est relativement irr\'eguli\`ere. Rappelons que $0.0018$ et $0.0050$ sont les valeurs num\'eriques de la constante $C_{\text{Tr}}^{(0)}$ apparaissant dans le corollaire \ref{rang0} pour les courbes $26\text{A}$ et $26\text{B}$.

\subsection{Courbes de rang analytique $1$}

\begin{figure}[hbp]
\begin{center}
\includegraphics[height=10cm]{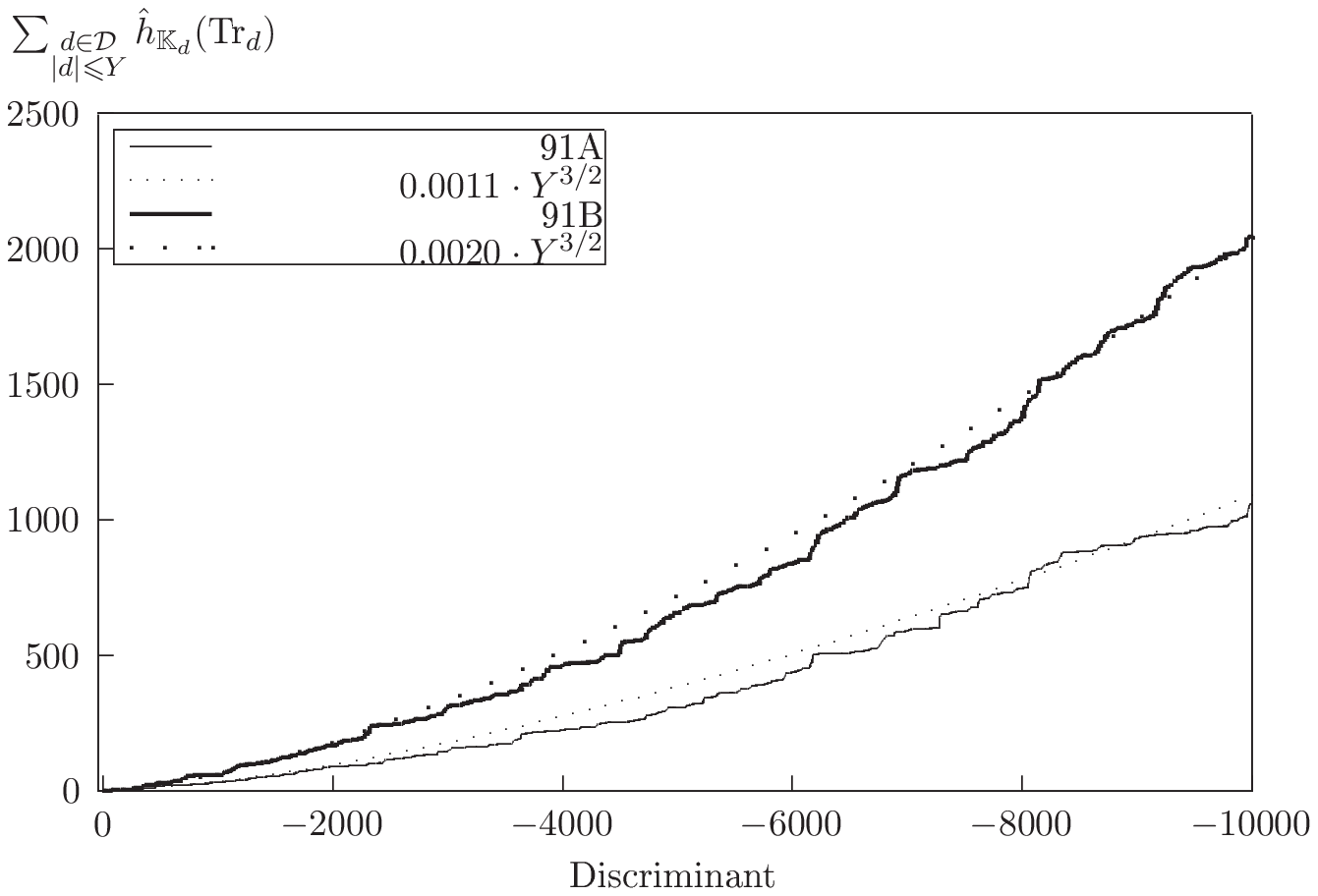}
\end{center}
\caption{Hauteur des traces en moyenne sur les courbes $91\text{A}$ et $91\text{B}$.}
\label{91AB-traces}
\end{figure}

\begin{figure}[htp]
\begin{center}
\includegraphics[height=10cm]{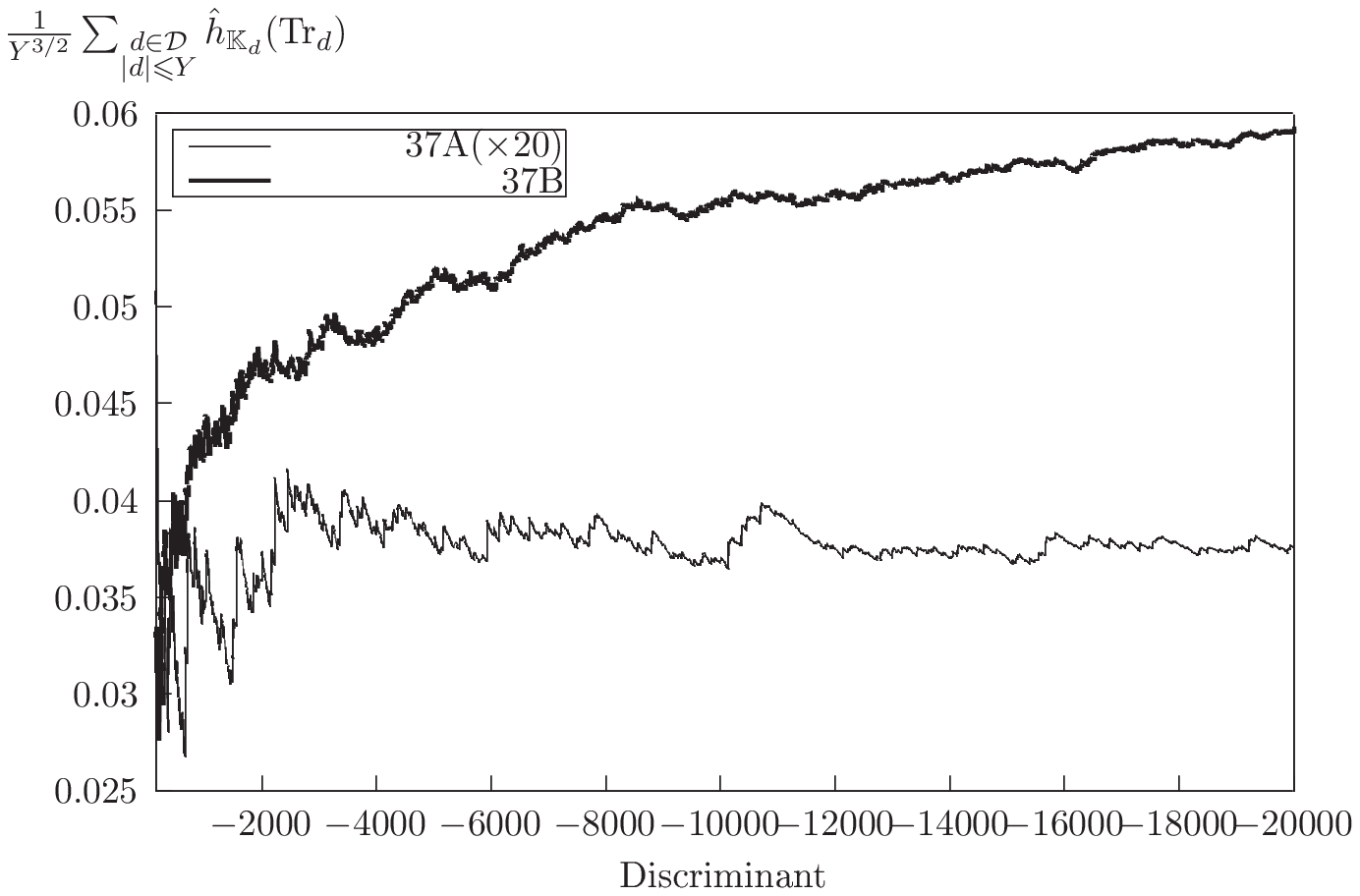}
\end{center}
\caption{Hauteur des traces en moyenne sur les courbes $37\text{A}$ et $37\text{B}$.}
\label{37AB-traces}
\end{figure}

Lorsque l'on consid\`ere une courbe de rang $1$, on est amen\'e \`a estimer la moyenne des valeurs en $1$ des fonctions $L$ tordues et non de leurs d\'eriv\'ees.
\begin{theorem}
\label{iwa1}
Si $E$ est une courbe elliptique rationnelle de conducteur $N$ sans facteurs carrés et de rang analytique $\,\,1$ et $F$ est une fonction lisse à support compact dans $\mathbb{R}_+$ et de moyenne strictement positive alors 
$$\sum_{d\in\mathcal{D}} L(E\vert\mathbb{Q}\times\chi_d,1)F\left(\frac{\vert d\vert}{Y}\right)=\alpha_NY+\mathcal{O}_\varepsilon\left(N^{\frac{23}{14}+\varepsilon}Y^{\frac{13}{14}+\varepsilon}\right)$$ 
pour tout $\varepsilon>0$ o\`u $\alpha_N$ est définie en \eqref{alphaN}.
\end{theorem}
\noindent{\textbf{Idée de preuve du théorème \ref{iwa1}.}} 
Il n'est pas difficile d'adapter la d\'emonstration de \cite{iwaniec} \`a ce cas. En reprenant les notations de l'article, il suffit de remplacer la fonction $V(X)$ d\'efinie au paragraphe 4 page $369$ par la fonction
$$\widetilde{V}(X)=e^{-X}.$$
On a alors\footnote{\`A noter une erreur de frappe dans \cite{iwaniec}, il s'agit bien de $\paf{X}{2\pi}^s$ et non de son inverse.},
$$\mathcal{A}\mn{(X,\chi_d)=\f{1}{2i\pi}\int_{(3/4)} L(s+1,E,\chi_d)\Gamma(s)\paf{X}{2\pi}^s ds}.$$
La majoration $\mathcal{A}\mn{(X,\chi_d)\ll  \sqrt{X}}$, qui d\'ecoule de l'in\'egalit\'e de Hölder et d'une estimation des $a_n$ tient toujours, et ainsi les majorations successives effectu\'ees dans la d\'emonstration ne posent pas de probl\`eme. La seule diff\'erence notable vient \`a la page 374 lors du calcul de $\mathcal{B}\mn{(X)}$. On a alors
$$res_{s=0} L(s+1)\Gamma(s)\paf{X}{2\pi}^s=L(1)$$
et il n'appara\^it pas de terme en $\log{X}$.
\begin{flushright}
$\blacksquare$
\end{flushright}
On en déduit comme pour le rang $0$ le corollaire suivant.
\begin{corollary}
\label{rang1}
Si $E$ est une courbe elliptique rationnelle de conducteur $N$ sans facteurs carrés et de rang analytique $1$ alors
\begin{equation*}
\sum_{\substack{d\in\mathcal{D} \\
\vert d\vert\leqslant Y}}\hat{h}_{\mathbb{K}_d}(\emph{Tr}_d)=C_{\emph{Tr}}^{(1)}Y^{3/2}+\mathcal{O}_\varepsilon\left(\frac{L(E\vert\mathbb{Q},1)}{\Omega_{E,N}}N^{\frac{23}{14}+\varepsilon}Y^{\frac{20}{14}+\varepsilon}\right)
\end{equation*}
pour tout $\varepsilon>0$ o\`u $C_{\emph{Tr}}^{(1)}$ est la constante d\'efinie par
$$C_{\emph{Tr}}^{(1)}:=\f{2}{\pi}c_N\mathcal{P}(1)L^\prime(E\vert\mathbb{Q},1)\f{L(\text{Sym}^2E,2)}{\pi\Omega_{E,N}}L_E$$
où $c_N$ est définie en \eqref{cN}, $L_E$ en \eqref{LE} et $\mathcal{P}$ en \eqref{P}.
\end{corollary}
\begin{remark}
\`A conducteur $N$ fixé, le terme principal dans \eqref{rang1} est $C_{\emph{Tr}}^{(1)}Y^{\frac{3}{2}}$ avec
\begin{equation*}
C_{\emph{Tr}}^{(1)}=\left(\frac{6}{\pi^3c_E(N)^2}\mathcal{P}(1)\prod_{\substack{p\in \mathcal{P} \\
p \mid 2N}}\left(1-\frac{1}{p^2}\right)^{-1}\right)\times L^\prime(E\vert\mathbb{Q},1)\,L_{E}\times\frac{\gamma(4N)}{N^2}\text{deg}(\Phi_{N,E})
\end{equation*}
et est donc proportionnel au degré de la paramétrisation modulaire.
\end{remark}
\begin{remark}
Il est intuitivement étonnant que les traces en moyennes des points de Heegner sur une courbe elliptique $E$ sont asymptotiquement plus grosses par un facteur logarithmique si le rang de la courbe elliptique est minimal.
\end{remark}
Les deux courbes elliptiques de conducteur $91$ ont m\^eme rang $1$ et on a repr\'esent\'e sur la figure \ref{91AB-traces} les hauteurs des traces Tr$_d$ sur ces deux courbes ainsi que la courbe th\'eorique donn\'ee par le corollaire ci-dessus. On constate que les courbes sont beaucoup plus irr\'eguli\`eres que dans le cas du rang $0$ (figure \ref{26A-traces}) m\^eme si elles suivent la courbe th\'eorique de tr\`es pr\`es.
La figure \ref{37AB-traces} illustre les diff\'erences de croissance des hauteurs des traces sur les courbes $37\text{A}$ (rang $1$) et $37\text{B}$ (rang $0$). On remarque que le comportement est tr\`es irr\'egulier. Ceci est en partie d\^u \`a la division par $Y^{3/2}$ qui rend les irr\'egularit\'es plus apparentes que dans la figure \ref{26A-traces}. De plus, pour la courbe $37\text{A}$ de rang $1$, il arrive fr\'equemment que la trace soit nulle ce qui <<casse la moyenne>>. Il est conjectur\'e que la proportion d'annulation de $L(E\vert\mathbb{Q}\times\chi_d,1)$ (ce qui correspond aux cas de trace nulle selon la formule de Gross-Zagier et la conjecture de Birch et Swinnerton-Dyer) tend vers $0$ lorsque le discriminant tend vers l'infini\footnote{Le type de symétrie de cette famille de fonctions $L$ est orthogonal impair.} mais cela ne se voit pas dans l'\'echelle de discriminants \'etudi\'ee.

\section{Estimation asymptotique de la hauteur des points de Heegner}

Nous d\'emontrons une formule asymptotique pour les hauteurs des points de Heegner P$_d$ semblable \`a celle que nous avons donn\'ee pour les traces. Soient $E$ une courbe elliptique définie sur $\mathbb{Q}$, de conducteur $N$ et $L(E\vert\Q,s):=\sum_{n\geqslant 1}a_nn^{-s}$ sa fonction $L$ (de rang analytique quelconque). On s'int\'eresse \`a la valeur moyenne des d\'eriv\'ees en $1$ des s\'eries de Dirichlet $L_d(E,s)$ d\'efinies par \eqref{lesd}.
\begin{theorem}
\label{traceana}
Si $E$ est une courbe elliptique rationnelle de conducteur $N$ sans facteurs carrés et de rang analytique quelconque et $F$ est une fonction lisse à support compact dans $\mathbb{R}_+$ et de moyenne strictement positive alors
$$\sum_{d\in\mathcal{D}}L^\prime_d(E,1)F\left(\frac{\vert d\vert}{Y}\right)=\widetilde{\alpha_N} Y\log{Y}+\widetilde{\beta_N}Y+\mathsf{Error}+\mathcal{O}_\varepsilon\left(N^{\frac{15}{4}+\varepsilon}Y^{\frac{19}{20}+\varepsilon}\right)$$
où
\begin{equation}
\label{conjerror}
\mathsf{Error}=\mathcal{O}_\varepsilon\left(NY\left(\log{(NY)}\right)^{\frac{1}{2}+\varepsilon}\right)
\end{equation}
pour tout $\varepsilon>0$ et o\`u
\begin{equation*}
\label{defbeta}
\widetilde{\alpha_N}:=c_N\widetilde{L}(1)\int_0^{+\infty}F(t)\mathrm{d}t\neq 0
\end{equation*}
et
\begin{equation*}
\widetilde{\beta_N}:=c_N\int_0^{+\infty}F(t)\left(\widetilde{L}^\prime(1)+\widetilde{L}(1)\left(\log{\left(\frac{Nt}{4\pi^2}\right)}-2\gamma\right)\right)\mathrm{d}t
\end{equation*}
où la constante $c_N$ est définie en \eqref{cN} avec
\begin{equation*}
\widetilde{L}(s):=\frac{L(\text{Sym}^2E,2s)}{\zeta^{(N)}(2s)}\widetilde{\mathcal{P}}(s)\times\begin{cases}
\frac{4}{3} & \text{ si $N$ est impair,} \\
1 & \text{ sinon}
\end{cases}
\end{equation*}
et
\begin{equation*}
\widetilde{\mathcal{P}}(s):=\prod_{\substack{p\in\mathcal{P} \\
p\nmid 2N}}\left(1+\left(1+\frac{1}{p}\right)^{-1}(p^{4s-2}-1)^{-1}\right).
\end{equation*}
\end{theorem}
\begin{remark}
\label{rq}
\`A conducteur $N$ fixé, il semble que l'on obtienne un développement asymptotique du premier moment par rapport à $Y$ à un seul terme et non à deux termes comme dans le Théorème \ref{iwa} de H. Iwaniec. Cependant, les résultats numériques décrits dans la dernière partie du paragraphe \ref{5} et notre intuition analytique nous permettent de conjecturer que
\begin{equation*}
\mathsf{Error}=o_\varepsilon(NY).
\end{equation*}
Pour pouvoir prouver cela, il faudrait notamment être en mesure de déterminer le comportement asymptotique de moyennes de la forme
\begin{equation*}
\sum_{\substack{1\leqslant u\leqslant U \\
1\leqslant v\leqslant V}}\sum_{d\in\mathcal{D}}a_{u,v}\chi_d(u)r_d(v)F\left(\frac{\vert d\vert}{Y}\right)
\end{equation*}
pour tous nombres réels strictement positifs $U$, $V$ et toute suite de nombres complexes $\left(a_{u,v}\right)_{\substack{1\leqslant u\leqslant U \\
1\leqslant v\leqslant V}}$ (se reporter également en page \pageref{fin}). Les auteurs projettent de s'intéresser dans un avenir proche à ce type de moyennes qui sont en réalité un cas particulier de quantités beaucoup plus générales. En outre, il ne fait aucun doute que la dépendance en $N$ dans $\mathcal{O}_\varepsilon\left(N^{\frac{15}{4}+\varepsilon}Y^{\frac{19}{20}+\varepsilon}\right)$ peut être améliorée en étant plus soigneux mais une croissance au plus polynomiale en le conducteur nous suffit.
\end{remark}
\noindent{\textbf{Preuve du théorème \ref{traceana}.}} B.H. Gross et D. Zagier (\cite{GZ}) ont prouvé que la s\'erie de Dirichlet $L_d(E,s)$ définie a priori sur $\Re{(s)}>\frac{3}{2}$ admet un prolongement holomorphe à $\mathbb{C}$ et satisfaisait l'\'equation fonctionnelle
$$\forall s\in\mathbb{C}, \quad \Lambda_d(E,s)=-\chi_d(N)\Lambda_d(E,2-s)$$
où $\Lambda_d(E,s):=(N\vert d\vert)^s\left((2\pi)^{-s}\Gamma(s)\right)^2L_d(E,s)$ est la série de Dirichlet complétée. Remarquons que comme $d$ est un carré modulo $N$, le signe de l'équation fonctionnelle vaut $-1$ d'où $L_d(E,1)=0$. Ceci va nous permettre d'exprimer $L^\prime_d(E,1)$ en terme de deux sommes convergeant exponentiellement vite en suivant une procédure analytique désormais classique (confer Théorème 5.3. de \cite{IwKo} pour plus de détails). Pour $X>0$, posons
$$V(X):=\f{1}{2i\pi}\int_{(3/4)} \Gamma(s)^2 X^{-s} ds,$$
et 
$$\mathcal{A}_d\mn{(E,X):=\f{1}{2i\pi}\int_{(3/4)}L_d(E,s+1)\Gamma(s)^2\paf{X}{4\pi^2}^s ds.}$$
Le développement de $L_d(E,s)$ en série de Dirichlet absolument convergente sur $\Re{(s)}>\frac{3}{2}$ assure que
$$\mathcal{A}_d\mn{(E,X)=\sum_{n=1}^{\infty} \f{a_n r_d(n)}{n}\sum_{(m,N)=1} \chi_d(m)\f{1}{m} V\paf{4\pi^2 nm^2}{X}.}$$
On d\'eplace la ligne d'int\'egration jusqu'à $\Re(s)=-3/4$ croisant un unique pole en $s=0$ de résidu \'egal \`a $L^\prime_d(E,1)$ puis on revient en $s=3/4$ par le changement de variables $s\mapsto -s$. L'\'equation fonctionnelle entraîne alors que
$$L^\prime_d(E,1)=\mathcal{A}_d\mn{(E,X)+}\mathcal{A}_d\mn{\left(E,\f{(Nd)^2}{X}\right),}$$
et en particulier que
$$L^\prime_d(E,1)=2\mathcal{A}_d\mn{(E,\vert d\vert N).}$$
Bornons $V(X)$ pour $X>0$ de la mani\`ere suivante :
\begin{align}
\label{V}
V(X)& =\f{1}{2i\pi}\int_{u=0}^{\infty}\f{e^{-u}}{u} \int_{v=0}^{\infty} e^{-v}\pa{ \int_{(3/4)} \paf{uv}{X}^s \f{1}{s} ds} dv du \nonumber \\
& \ll \int_{0}^{\infty} \f{e^{-(u+X/u)}}{u}du \nonumber \\
& \ll X^{-1/4}\exp{\left(-2\sqrt{X}\right)}
\end{align}
et en fait $X^jV^{(j)}(X)\ll_jX^{-1/4}\exp{\left(-2\sqrt{X}\right)}$ pour tout entier naturel $j$. Posons
$$S_N(Y):=\sum_{d\in\mathcal{D}}L^\prime_d(E,1)F\left(\frac{\vert d\vert}{Y}\right).$$
Comme $d$ est dans $\mathcal{D}$, $d$ est un carr\'e modulo $4$ et est premier \`a $4$ donc $d$ est congru à $1$ modulo $4$. Ainsi, $\mathcal{O}_d=\mathbb{Z}+\frac{1+\sqrt{d}}{2}\mathbb{Z}$ et un calcul élémentaire montre que
\begin{equation*}
r_d(n)=\#\left\{(u,v)\in\left(\N^{*}\times\Z\right)\cup\left(\{0\}\times\N\right), u^2+\vert d\vert v^2=4n\right\}.
\end{equation*}
On observe que si $n$ est un carr\'e alors $(2\sqrt{n},0)$ est une solution de l'équation ci-dessus alors qu'il n'existe pas de solutions de la forme $(0,*)$ ce qui prouve que
\begin{equation*}
r_d(n)=1+\#\left\{v\in\Z^*, \vert d\vert v^2=4n\right\}+\#\left\{(u,v)\in\left(\N^*\times\Z^*\right), u^2+\vert d\vert v^2=4n\right\}:=1+r_d^\prime(n)
\end{equation*}
et que
$$\sum_{\substack{d\in \mathcal{D} \\
Y\ll\vert d\vert\ll Y}}r_d(n)\geqslant\#\left\{d\in\mathcal{D}, Y\ll\vert d\vert \ll Y\right\}\sim_{Y\rightarrow +\infty}CY$$
pour une constante absolue $C>0$. Si $n$ n'est pas un carr\'e alors $r_d(n)=r_d^\prime(n)$. Dans chacun des cas, si $(u,v)$ est une solution contribuant à $r_d^\prime(n)$ pour un $d$ inférieur à $Y$ dans $\mathcal{D}$ alors $u\leqslant 2\sqrt{n}$ et \`a chaque tel $u$ correspond au plus deux couples $(d,v)$ (car $d$ est suppos\'e sans facteurs carr\'es) d'o\`u
\begin{equation}
\label{rdprime}
\sum_{\substack{d\in\mathcal{D} \\
Y\ll\vert d\vert\ll Y}}r_d^\prime(n)\leqslant 4\sqrt{n}.
\end{equation}
On \'ecrit alors\footnote{Au cours de cette preuve, $\mathsf{TP}_i$ désignera la quantité d'où provient la contribution principale et $\mathsf{Err}_i$ un terme d'erreur.} $S_N(Y):=\mathsf{TP}_1+\mathsf{Err}_1$ o\`u
$$\mathsf{TP}_1:=2\sum_{n\geqslant 1}\f{a_{n^2}}{n^2}\sum_{(m,N)=1} \f{1}{m} \sum_{d\in \mathcal{D}}\chi_d(m)V\paf{4\pi^2 n^2 m^2}{N\vert d\vert}F\left(\frac{\vert d\vert}{Y}\right)$$
et
\begin{equation}
\label{err1}
\mathsf{Err}_1:=2\sum_{d\in \mathcal{D}}\sum_{n\geqslant 1}\f{a_{n}r'_d(n)}{n}\sum_{(m,N)=1}\f{1}{m}\chi_d(m)V\paf{4\pi^2 n m^2}{N\vert d\vert}F\left(\frac{\vert d\vert}{Y}\right).
\end{equation}
\noindent{\textbf{Estimation du terme d'erreur $\mathbf{\mathsf{Err}_1}$.}}
Découpons $\mathsf{Err}_1$ de la façon suivante:
\begin{equation*}
\mathsf{Err}_1=2\sum_{(m,N)=1}\sum_{1\leqslant n\leqslant\frac{NY\psi(NY)}{m^2}}\cdots+2\sum_{(m,N)=1}\sum_{n>\frac{NY\psi(NY)}{m^2}}\cdots:=\mathsf{Error}+\mathsf{Err}_3
\end{equation*}
pour toute fonction $\psi$ positive et tendant vers $0$ en $+\infty$. Les estimations \eqref{V} et \eqref{rdprime} assurent que
\begin{equation*}
\mathsf{Error}\ll(NY)^{1/4}\sum_{1\leqslant m\leqslant\sqrt{NY\psi(NY)}}\f{1}{m^{3/2}}\sum_{1\leqslant n\leqslant\frac{NY\psi(NY)}{m^2}}\frac{\vert a_n\vert}{n^{\frac{3}{4}}}\exp{\left(-4\pi\frac{m}{\sqrt{NY}}\sqrt{n}\right)}.
\end{equation*}
L'inégalité de Cauchy-Schwarz entraîne que le carré de la somme en $n$ est borné par
\begin{equation*}
\left(\sum_{1\leqslant n\leqslant\frac{NY\psi(NY)}{m^2}}\f{a_n^2}{n^2}\right)\left(\sum_{1\leqslant n\leqslant\frac{NY\psi(NY)}{m^2}}\sqrt{n}\exp{\left(-8\pi\frac{m}{\sqrt{NY}}\sqrt{n}\right)}\right).
\end{equation*}
La première somme est estimée par $\mathcal{O}\left(\log{(NY\psi(NY))}\right)$ alors que la deuxième est trivialement inférieure à
\begin{equation*}
\int_{1}^{\frac{NY\psi(NY)}{m^2}}\sqrt{t}\exp{\left(-8\pi\frac{m}{\sqrt{NY}}\sqrt{t}\right)}\mathrm{d}t\ll\left(\frac{NY\psi(NY)}{m^2}\right)^{\frac{3}{2}}.
\end{equation*}
Ainsi, on a prouvé que $\mathsf{Error}\ll NY(\psi(NY))^{\frac{3}{4}}\left(\log{(NY\psi(NY))}\right)^{\frac{1}{2}}$. De la m\^eme mani\`ere, 
\begin{equation*}
\mathsf{Err}_3\ll(NY)^{1/4}\sum_{m\geqslant 1}\f{1}{m^{3/2}}\left(\sum_{n>\frac{NY\psi(NY)}{m^2}}\frac{a_n^2}{n^{2+\varepsilon}}\right)^{\frac{1}{2}}\left(\sum_{n>\frac{NY\psi(NY)}{m^2}}n^{\frac{1}{2}+\varepsilon}\exp{\left(-8\pi\frac{m}{\sqrt{NY}}\sqrt{n}\right)}\right)^{\frac{1}{2}}
\end{equation*}
pour tout $\varepsilon>0$. Une intégration par parties assure que
\begin{equation*}
\mathsf{Err}_3\ll(NY)^{1+\varepsilon}(\psi(NY))^{\frac{1}{2}+\varepsilon}\exp{\left(-\frac{1}{2}\sqrt{\psi(NY)}\right)}
\end{equation*}
et on choisit alors $\psi(x):=(\log{x})^a$ avec $2\varepsilon<a<\frac{2}{3}$ de sorte que $\mathsf{Err}_3=o(\mathsf{Error})$ et que $\mathsf{Error}\ll (NY)(\log{(NY)})^{\frac{3a}{4}+\frac{1}{2}}=o((NY)\log{(NY)})$.\newline
\noindent{\textbf{Contribution du terme principal $\mathbf{\mathsf{TP}_1}$.}}
Déterminons le comportement asymptotique de $\mathsf{TP}_1$ en appliquant une méthode développée dans \cite{iwaniec}:
\begin{itemize}
\item
la condition $d$ sans facteurs carrés est supprimée en introduisant $\sum_{a^2\mid d}\mu(a)$ 
puis la somme est coupée selon la taille des diviseurs $a$ de $d$ ($a\leqslant A$ et $a>A$) sachant que l'on revient à des discriminants sans facteurs carrés dans le cas des grands diviseurs;
\item
pour tout entier $m=m_1m_2^2$ avec $(m_1m_2,N)=1$ et $m_1$ sans facteurs carrés, remarquons que $\chi_d(m)=\chi_d(m_1)$ si $(m_2,d)=1$ (et $0$ sinon) puis que le développement de Fourier du caractère $\chi_.(m_1)$ en terme de caractères additifs de module $m_1$ s'écrit
\begin{equation*}
\chi_d(m_1)=\frac{\overline{\varepsilon_{m_1}}}{\sqrt{m_1}}\sum_{0\leqslant\vert r\vert<\frac{m_1}{2}}\chi_{Nr}(m_1)e\left(\frac{\overline{N}rd}{m_1}\right)
\end{equation*}
où $\overline{N}$ est l'inverse de $N$ modulo $m_1$ et $\varepsilon_{m_1}$ est le signe de la somme de Gauss de $\chi_.(m_1)$.
\end{itemize}
La contribution principale provient alors du terme $r=0$ pour lequel $\chi_0(m_1)$ vaut $0$ si $m_1>1$ et $1$ sinon. En résumé,
\begin{equation*}
\mathsf{TP}_1=\mathsf{TP}_2+\mathsf{Err}_4+\mathsf{Err}_5
\end{equation*}
où
\begin{eqnarray*}
\mathsf{TP}_2 & :=  & 2\sum_{n\leqslant 1}\f{a_{n^2}}{n^2}\sum_{\substack{a\leqslant A \\
(a,4N)=1}}\mu(a)\sum_{(m_2,aN)=1}\f{1}{m_2^2}\sum_{q\mid m_2}\mu(q)\sum_{\substack{qd\in\mathcal{D}^\prime \\
(d,m_2)=1}}V\paf{4\pi^2 n^2 m_2^4}{Na^2\vert d\vert q}F\left(\frac{a^2\vert d\vert q}{Y}\right), \\
\mathsf{Err}_5 & := & 2\sum_{n\leqslant 1}\f{a_{n^2}}{n^2}\sum_{(b,4N)=1}\sum_{\substack{a\mid b \\
a>A}}\mu(a)\sum_{(m,N)=1} \f{1}{m} \sum_{d\in \mathcal{D}}\chi_{b^2d}(m)V\paf{4\pi^2 n^2 m^2}{Nb^2\vert d\vert}F\left(\frac{b^2\vert d\vert}{Y}\right)
\end{eqnarray*}
et
\begin{multline*}
\mathsf{Err}_4:=2\sum_{n\leqslant 1}\f{a_{n^2}}{n^2}\sum_{\substack{a\leqslant A \\
(a,4N)=1}}\mu(a)\sum_{\substack{m=m_1m_2^2 \\
(m,aN)=1}}\f{\mu^2(m_1)}{m}\sum_{q\mid m_2}\mu(q)\sum_{\substack{qd\in\mathcal{D}^\prime \\
(d,m_2)=1}} \\
\frac{\overline{\varepsilon_{m_1}}}{\sqrt{m_1}}\sum_{1\leqslant\vert r\vert<\frac{m_1}{2}}\chi_{Nrq}(m_1)e\left(\frac{\overline{N}rd}{m_1}\right)V\paf{4\pi^2 n^2 m^2}{Na^2\vert d\vert q}F\left(\frac{a^2\vert d\vert q}{Y}\right)
\end{multline*}
avec
\begin{equation*}
\mathcal{D}^\prime:= \mn{\{d\in\mathbb{Z}_-^*, d\equiv \nu^2 \mod 4N, (\nu,4N)=1\}}.
\end{equation*}
\noindent{\textbf{Estimation du terme d'erreur $\mathbf{\mathsf{Err}_5}$.}}
Pour commencer, l'inégalité de Hölder implique que
\begin{equation*}
\mathsf{Err}_5\ll\sum_{n\ll\left(NY\right)^{\frac{1}{2}}}\f{\left\vert a_{n^2}\right\vert}{n^2}\sum_{\substack{b\geqslant 1 \\
a\mid b \\
a>A}}\left(\sum_{\substack{d\in \mathcal{D} \\
\vert d\vert\ll\frac{Y}{b^2}}}1\right)^{\frac{3}{4}}\left(\sum_{\substack{d\in \mathcal{D} \\
\vert d\vert\ll\frac{Y}{b^2}}}\chi_{b^2d}(m)\left\vert\sum_{m\ll\left(\frac{NY}{n^2}\right)^{\frac{1}{2}}} \f{1}{m}\right\vert^4\right)^{\frac{1}{4}}
\end{equation*}
d'où trivialement
\begin{equation*}
\mathsf{Err}_5\ll Y^{\frac{3}{4}}\sum_{n\ll\left(NY\right)^{\frac{1}{2}}}\f{\left\vert a_{n^2}\right\vert}{n^2}\sum_{\substack{b\geqslant 1 \\
a\mid b \\
a>A}}\frac{1}{b^{\frac{3}{2}}}\left(\sum_{\substack{d\in \mathcal{D} \\
\vert d\vert\ll\frac{Y}{b^2}}}\chi_{b^2d}(m)\left\vert\sum_{m\ll\left(\frac{NY}{n^2}\right)^{\frac{1}{2}}} \f{1}{m}\right\vert^2\right)^{\frac{1}{2}}.
\end{equation*}
L'inégalité du grand crible (\cite{Bo}) pour les caractères réels assure alors que
\begin{equation*}
\mathsf{Err}_5\ll_\varepsilon \frac{Y^{\frac{5}{4}+\varepsilon}}{A^{\frac{3}{2}}}+\frac{N^{\frac{1}{4}+\varepsilon}Y^{1+\varepsilon}}{A^{\frac{1}{2}}}
\end{equation*}
pour tout $\varepsilon>0$.\newline
\noindent{\textbf{Estimation du terme d'erreur $\mathbf{\mathsf{Err}_4}$.}}
Posons $\Delta:=\inf{\left(\frac{1}{2},\frac{a^2q}{Y^{1-\varepsilon}}\right)}$ pour tout nombre réel $\varepsilon>0$ et découpons $\mathbf{\mathsf{Err}_4}$ selon que la sommation en $r$ est restreinte par
\begin{eqnarray*}
1\leqslant\vert r\vert<\Delta m_1 & \rightsquigarrow & \mathsf{Err}_6, \\
\Delta m_1\leqslant\vert r\vert<\frac{m_1}{2} & \rightsquigarrow & \mathsf{Err}_7.
\end{eqnarray*}
On estime $\mathsf{Err}_7$ en bornant la somme sur les discriminants grâce au lemme 2 page 372 de \cite{iwaniec} puis trivialement la somme en $n$ et $m$ ce qui entraine que
\begin{equation*}
\mathsf{Err}_7\ll_\varepsilon\gamma(4N)Y^\varepsilon N^{\frac{5}{4}+\varepsilon}\frac{ \inf{(A,Y^{\frac{1-\varepsilon}{2}}) }}{Y^{\varepsilon-\frac{1}{4}}}+\gamma(4N)Y^\varepsilon N^{\frac{5}{4}+\varepsilon}\frac{A^3}{Y^{\frac{3}{4}}}
\end{equation*}
où $\gamma(4N)$ est le cardinal de l'ensemble des classes d'équivalence de $\mathcal{D}^\prime$ modulo $4N$. On estime $\mathsf{Err}_6$ de façon triviale par
\begin{equation*}
\mathsf{Err}_6\ll_\varepsilon N^{\frac{3}{8}+\varepsilon}Y^{\frac{1}{2}+\varepsilon}A^{\frac{3}{2}}.
\end{equation*}
\noindent{\textbf{Contribution du terme principal $\mathbf{\mathsf{TP}_2}$.}}
Intéressons-nous au terme principal $\mathsf{TP}_2$ et plus précisément à la somme sur les discriminants intervenant dans cette somme. Pour cela, on note $\mathcal{D}^\prime(4N)$ l'ensemble des classes d'équivalence de $\mathcal{D}^\prime$ modulo $4N$ et on se souvient que $\#\mathcal{D}^\prime(4N)=\gamma(4N)$. La formule de Poisson assure que
\begin{multline*}
\sum_{qd\in\mathcal{D}^\prime}V\paf{4\pi^2 n^2 m_2^4}{Na^2\vert d\vert q}F\left(\frac{a^2\vert d\vert q}{Y}\right)=\frac{Y}{4Na^2q}\sum_{[d_0]\in\mathcal{D}^\prime(4N)}\sum_{\ell\in\mathbb{Z}} \\
\times\int_{\mathbb{R}}V\paf{4\pi^2 n^2 m_2^4}{NY\left\vert \frac{a^2qd_0}{Y}+t\right\vert}F\left(\left\vert\frac{a^2qd_0}{Y}+t\right\vert\right)e\left(\frac{Y\ell}{4Na^2q}t\right)\mathrm{d}t.
\end{multline*}
On isole alors le terme $\ell=0$ et on effectue deux intégrations par parties pour chaque terme $\ell\neq 0$ afin de rendre absolument convergente la série en $\ell$ (il ne reste pas de termes entre crochets car $F$ est à support compact). On obtient alors
\begin{multline*}
\mathsf{TP}_2=\frac{\gamma(4N)Y}{2N}\int_0^{+\infty}F(t)\left(\sum_{n\geqslant 1}\f{a_{n^2}}{n^2}\sum_{\substack{a\leqslant A \\
(a,4N)=1}}\frac{\mu(a)}{a^2}\sum_{(m_2,aN)=1}\f{1}{m_2^2}\sum_{q\mid m_2}\frac{\mu(q)}{q}V\paf{4\pi^2 n^2 m_2^4}{NYt}\right)\mathrm{d}t \\
+\mathcal{O}\left(\frac{N^{\frac{5}{4}}\gamma(4N)A^3}{Y^{\frac{3}{4}}}\right)
\end{multline*}
Finalement,
$$\mathsf{TP}_2=c_N Y \int_0^{+\infty}F(t)\mathcal{B}(NYt)dt+\mathcal{O}\left(\frac{N^{\frac{1}{4}}\gamma(4N)Y^{\frac{5}{4}}}{A}+\frac{N^{\frac{5}{4}}\gamma(4N)A^3}{Y^{\frac{3}{4}}}\right)$$
avec
$$\mathcal{B}\mn{(X)=\sum_{n\geqslant 1}\f{a_{n^2}}{n^2}}\sum_{(m,N)=1}\f{b_{m^2}}{m^2}V\paf{4\pi^2n^2m^4}{X}$$
et
$$b_m=\prod_{\substack{p\in\mathcal{P} \\
p\mid m \\
p\neq 2}} \pa{1+\f{1}{p}}^{-1}.$$
En revenant à la définition intégrale de la fonction $V$, on remarque que
\begin{equation*}
\mathcal{B}(X)=\f{1}{2i\pi}\int_{(3/4)} \Gamma(s)^2 X^{-s}L(s+1)\mathrm{d}s
\end{equation*}
et que le produit Eulérien intervenant dans la fonction $L$ est absolument convergeant sur $\Re{(s)}>\frac{3}{4}$ et y définit une fonction holomorphe. En décalant le contour jusqu'à $\left(-\frac{1}{4}+\varepsilon\right)$ pour tout $\varepsilon>0$, on ne croise qu'un pôle en $s=0$ ce qui prouve que
\begin{equation*}
\mathcal{B}(X)=-2(\gamma+\log{(2\pi)})L(1)+L^\prime(1)+L(1)\log{(X)}+\mathcal{O}_\varepsilon\left(\left(\frac{N}{X}\right)^{\frac{1}{4}+\varepsilon}\right).
\end{equation*}
\noindent{\textbf{Bilan et choix des paramètres.}}
On a prouvé que
\begin{equation*}
S_N(Y)=\widetilde{\alpha_N} Y\log{Y}+\widetilde{\beta_N}Y+\mathcal{O}_\varepsilon\left(NY\left(\log{(NY)}\right)^{\frac{1}{2}+\varepsilon}\right)+\mathsf{Err}
\end{equation*}
où
\begin{equation*}
\mathsf{Err}\ll_\varepsilon(NY)^\varepsilon\left(\frac{Y^{\frac{5}{4}}}{A^{\frac{3}{2}}}+\frac{N^{\frac{1}{4}}Y}{A^{\frac{1}{2}}}+\frac{N^{\frac{9}{4}}\inf{(A,Y^{\frac{1-\varepsilon}{2}})}}{Y^{\varepsilon-\frac{1}{4}}}+\frac{N^{\frac{9}{4}}A^3}{Y^{\frac{3}{4}}}+N^{\frac{3}{8}}Y^{\frac{1}{2}+\varepsilon}A^{\frac{3}{2}}+\frac{N^{\frac{5}{4}}Y^{\frac{5}{4}}}{A}\right)
\end{equation*}
et on choisit alors $A:=N^{\frac{1}{2}}Y^{\frac{3}{10}-\frac{2\varepsilon}{5}}$ ce qui achève la preuve.
\begin{flushright}
$\blacksquare$
\end{flushright}
Appliquons finalement la formule de Gross-Zagier \eqref{gzpoint} pour obtenir une estimation asymptotique de la hauteur en moyenne des points de Heegner de la m\^eme forme que celle que l'on avait obtenue pour les traces.
\begin{corollary}
\label{hpointstheo}
Si $E$ est une courbe elliptique rationnelle de conducteur $N$ sans facteurs carr\'es et de rang analytique quelconque alors
$$\sum_{\substack{d\in\mathcal{D} \\
\vert d\vert\leqslant Y}}\hat{h}_{\mathbb{H}_d}(\emph{P}_d)=C_{\emph{P}}Y^{\frac{3}{2}}\log{Y}+C_{\emph{P}}^\prime Y^{\frac{3}{2}}+\frac{1}{3\Omega_{E,N}}\sqrt{Y}\mathsf{Error}+\mathcal{O}_\varepsilon\left(N^{\frac{15}{4}+\varepsilon}Y^{\frac{29}{20}+\varepsilon}\right)$$
où
\begin{equation*}
\frac{1}{3\Omega_{E,N}}\sqrt{Y}\mathsf{Error}=\mathcal{O}_\varepsilon\left(NY^{\frac{3}{2}}\left(\log{(NY)}\right)^{\frac{1}{2}+\varepsilon}\right)
\end{equation*}
pour tout $\varepsilon>0$ et o\`u $C_{\emph{P}}$ est la constante d\'efinie par
$$C_{\emph{P}}:=\f{2}{\pi}c_N\mathcal{Q}(N)\f{L(\text{Sym}^2 E,2)}{\pi\Omega_{E,N}}\prod_{\substack{p\in\mathcal{P} \\
p\mid N}}\left(1-\frac{1}{p^2}\right)^{-1}$$
avec
\begin{equation*}
\mathcal{Q}(N):=\prod_{\substack{p\in\mathcal{P} \\
p\nmid 2N}}\left(1+\left(1+\frac{1}{p}\right)^{-1}(p^{2}-1)^{-1}\right)\times\begin{cases}
\frac{4}{3} & \text{ si $N$ est impair,} \\
1 & \text{ sinon}
\end{cases}
\end{equation*}
et
\begin{equation*}
C_{\emph{P}}^\prime:=C_{\emph{P}}\left(\log{\left(\frac{N}{4\pi^2}\right)}-\frac{2}{3}-2\gamma\right)+\frac{c_N}{3\Omega_{E,N}}\widetilde{L}^\prime(1).
\end{equation*}
\end{corollary}
\begin{remark}
En accord avec la remarque \ref{rq}, on peut conjecturer que
\begin{equation*}
\frac{1}{3\Omega_{E,N}}\sqrt{Y}\mathsf{Error}=o_\varepsilon\left(NY^{\frac{3}{2}}\right)
\end{equation*}
et le corollaire précédent semble alors nous munir d'un développement asymptotique à deux termes de la hauteur en moyenne des points de Heegner.
\end{remark}
\begin{remark}
En rempla\c{c}ant $c_N$ et $L(\text{Sym}^2 E,2)$ par leur expression on peut r\'e\'ecrire 
$$C_{\text{P}}=\left(\paf{8}{\pi^3c_E(N)^2}\mathcal{Q}\mn{(N)}\prod_{\substack{p\in\mathcal{P} \\
p\mid N}}\pa{1-\f{1}{p^2}}^{-2}\right)\f{\gamma(4N)}{N^2}\text{deg}(\phi_{N,E}).$$
Ainsi, contrairement \`a la constante $C_{\text{Tr}}$ intervenant lorsque l'on consid\`ere les traces, \`a conducteur fix\'e $C_{\text{P}}$ ne d\'epend que du degr\'e de la param\'etrisation modulaire, puisque le produit $\mathcal{Q}\mn{(N)}$ ne d\'epend que de $N$. Par contre, lorsque l'on varie le conducteur, il n'y a plus une d\'ependance directe sur le degr\'e. Il est clair que $\mathcal{Q}\mn{(N)\geqslant 1}$ ; d'autre part
$$\mathcal{Q}\mn{(N) < \f{4}{3} \prod_p \pa{1+\f{1}{p^2}} < \f{4}{3} \zeta(2)}.$$
Le produit Eulerien $\mathcal{Q}\mn{(N)}$ est donc compris entre $1$ et $2$, et il joue un relativement faible r\^ole dans l'expression de $C_{\text{P}}$. On a ainsi
$$ 1 \leqslant \prod_{p|N}\pa{1-\f{1}{p^2}}^{-2}\mathcal{Q}\mn{(N) \leqslant \f{4}{3}\zeta(2)^3 < 6.}$$
Le terme principal, du moins si l'on s'int\'eresse \`a des valeurs asymptotiques du conducteur ou du degr\'e de la param\'etrisation modulaire, est donc $\f{\gamma(4N)}{N^2}\text{deg}(\phi_{N,E})$. Selon la conjecture du degr\'e (cf. \cite{mur2} et \cite{De} page 35) qui est équivalente à une des formes de la conjecture abc, on aurait $\text{deg}(\Phi_{N,E})\ll_\varepsilon N^{2+\epsilon}$ pour tout $\epsilon>0$. Comme $\gamma(4N)\ll N$, cela donne une borne sup\'erieure sur la croissance des hauteurs des points P$_d$ lorsque $N$ tend vers $+\infty$ avec $Y$. On sait d'autre part qu'il existe des familles de courbes de $j$-invariant borné (\cite{De} page 50) pour lesquelles $\text{deg}(\Phi_{N,E})\gg N^{\frac{7}{6}}\log{N}$ ce qui donne une borne inf\'erieure sur la vitesse de croissance des hauteurs lorsque $N$ tend vers $+\infty$ avec $Y$.
\end{remark}
\begin{remark}
Remarquons finalement que, même à conducteur fixé, la constante $C_{\emph{P}}^\prime$ dépend de la courbe elliptique $E$ et pas seulement du degré de la paramétrisation modulaire de $E$. Par contre, il ne semble pas être possible d'obtenir une estimation satisfaisante de la taille de $C_{\emph{P}}^\prime$ par rapport au niveau.
\end{remark}

\section{Analyse des résultats théoriques et numériques}
\label{5}
Après avoir donné quelques valeurs numériques des constantes en jeu, nous donnons des résultats expérimentaux illustrant les formules théoriques.

\subsection{Quelques valeurs num\'eriques}
\label{valnum}

\subsubsection{Valeurs numériques de $C_{\text{Tr}}$ et de $C_{\text{P}}$}
Le tableau \ref{constantes} regroupe les valeurs des trois constantes $C_{\text{Tr}}^{(0)}$, $C_{\text{Tr}}^{(1)}$ et $C_{\text{P}}$ r\'egissant le comportement en moyenne des hauteurs des points de Heegner et de leurs traces (selon les corollaires \ref{rang0}, \ref{rang1} et \ref{hpointstheo}) pour toutes les courbes elliptiques de conducteur sans facteurs carr\'es et inf\'erieur \`a $100$. Les valeurs de $C_{\text{Tr}}$ et de $C_{\text{P}}$ ont \'et\'e multipli\'ees par $10^3$ pour une meilleure lisibilit\'e.
\begin{table}[htbp]
\begin{center}
$\begin{array}{ccc}
\begin{array}{|c|c|c|c|c|}
\hline
\text{Courbe}&\text{Rang}& C_{\text{Tr}}\times 10^3 & C_{\text{P}}\times 10^3 & C_{\text{P}}/C_{\text{Tr}} \\
\hline
11-1 & 0 & 3.33 & 17.0 & 5.11 \\
14-1 & 0 & 1.07 & 6.39 & 5.92 \\
15-1 & 0 & 0.904 & 4.40 & 4.87 \\
17-1 & 0 & 3.21 & 11.3 & 3.52 \\
19-1 & 0 & 3.10 & 10.2 & 3.29 \\
21-1 & 0 & 1.21 & 3.28 & 2.71 \\
26-1 & 0 & 1.80 & 7.29 & 4.03 \\
26-2 & 0 & 5.00 & 7.29 & 1.45 \\
30-1 & 0 & 0.621 & 2.20 & 3.54 \\
33-1 & 0 & 2.48 & 6.56 & 2.64 \\
34-1 & 0 & 5.11 & 5.67 & 1.10 \\
35-1 & 0 & 1.95 & 4.29 & 2.19 \\
37-1 & 1 & 1.86 & 10.7 & 5.75 \\
37-2 & 0 & 5.64 & 10.7 & 1.90 \\
38-1 & 0 & 4.70 & 15.3 & 3.25 \\
38-2 & 0 & 5.09 & 5.10 & 1.00 \\
39-1 & 0 & 1.54 & 3.75 & 2.43 \\
42-1 & 0 & 2.30 & 3.28 & 1.42 \\
43-1 & 1 & 1.90 & 9.27 & 4.87 \\
46-1 & 0 & 3.16 & 10.6 & 3.36 \\
51-1 & 0 & 2.12 & 2.91 & 1.37 \\
53-1 & 1 & 1.87 & 7.55 & 4.02 \\
55-1 & 0 & 2.33 & 2.85 & 1.22 \\
57-1 & 1 & 0.923 & 5.24 & 5.68 \\
57-2 & 0 & 3.83 & 3.93 & 1.02 \\
57-3 & 0 & 8.24 & 15.7 & 1.90 \\

\hline
\end{array}
& &
\begin{array}{|c|c|c|c|c|}
\hline
\text{Courbe}&\text{Rang}& C_{\text{Tr}} & C_{\text{P}} & C_{\text{P}}/C_{\text{Tr}} \\
\hline
58-1 & 1 & 1.14 & 6.80 & 5.91 \\
58-2 & 0 & 9.47 & 6.80 & 0.718 \\
61-1 & 1 & 1.97 & 6.58 & 3.33 \\
62-1 & 0 & 4.85 & 3.18 & 0.657 \\
65-1 & 1 & 0.593 & 2.44 & 4.12 \\
66-1 & 0 & 1.01 & 2.18 & 2.15 \\
66-2 & 0 & 2.01 & 2.18 & 1.08 \\
66-3 & 0 & 25.2 & 10.9 & 0.433 \\
67-1 & 0 & 12.4 & 15.0 & 1.20 \\
69-1 & 0 & 2.40 & 2.18 & 0.909 \\
70-1 & 0 & 2.49 & 2.14 & 0.861 \\
73-1 & 0 & 6.90 & 8.27 & 1.19 \\
77-1 & 1 & 1.24 & 4.26 & 3.42 \\
77-2 & 0 & 15.4 & 21.3 & 1.37 \\
77-3 & 0 & 5.20 & 6.39 & 1.23 \\
78-1 & 0 & 4.47 & 18.7 & 4.18 \\
79-1 & 1 & 1.97 & 5.10 & 2.58 \\
82-1 & 1 & 1.15 & 4.85 & 4.19 \\
83-1 & 1 & 1.93 & 4.85 & 2.51 \\
85-1 & 0 & 3.00 & 3.80 & 1.26 \\
89-1 & 1 & 1.90 & 4.53 & 2.37 \\
89-2 & 0 & 10.5 & 11.3 & 1.07 \\
91-1 & 1 & 1.09 & 3.65 & 3.34 \\
91-2 & 1 & 2.03 & 3.65 & 1.79 \\
94-1 & 0 & 4.11 & 2.12 & 0.516 \\
     &   &      &      &       \\
\hline
\end{array}\\
\end{array}$
\end{center}
\caption{Valeurs numériques des constantes $C_{\text{Tr}}^{(0)}$, $C_{\text{Tr}}^{(1)}$ et $C_{\text{P}}$.}
\label{constantes}
\end{table}
\normalsize
Le rapport entre la constante gouvernant le comportement des hauteurs des points et celle donnant celui des traces donn\'e dans la derni\`ere colonne n'a que peu de sens dans le cas d'une courbe elliptique de rang $1$ car les points y sont asymptotiquement <<plus gros>> que les traces d'un facteur $\log{Y}$ selon les corollaires \ref{rang1} et \ref{hpointstheo}. Il est int\'eressant de voir que cette constante prend \`a la fois des valeurs plus grandes (<<les points sont plus gros>>) et plus petites (<<les traces sont plus grosses>>) que $1$.
\begin{remark} 
Si l'on poursuit le calcul sur les 200 premi\`eres courbes elliptiques de conducteur sans facteurs carr\'es alors on obtient un rapport moyen de $1.5$ environ et ce rapport tend \`a d\'ecro\^itre. Il n'y a donc pas de raison a priori de croire qu'il soit plus souvent plus grand ou petit que $1$.
\end{remark} 

\subsubsection{\'Etude plus fine du rapport $C_{\text{P}}/C_{\text{Tr}}^{(0)}$}
\begin{figure}[htbp]
\begin{center}
\includegraphics[height=10cm]{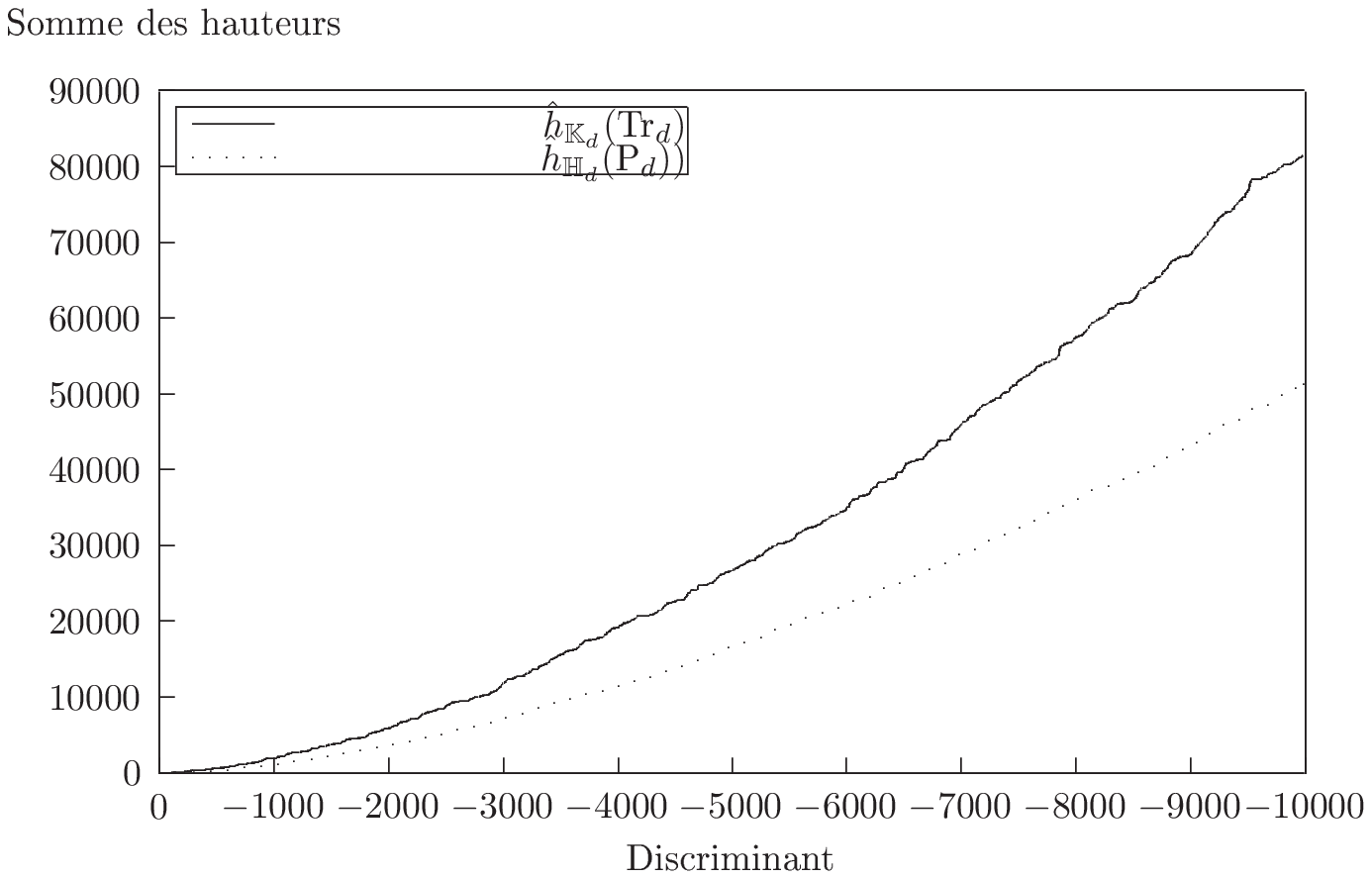}
\end{center}
\caption{Somme des hauteurs des traces et des points sur la courbe $58\text{B}$.}
\label{58B-points-traces}
\end{figure}

\begin{figure}[htbp]
\begin{center}
\includegraphics[height=10cm]{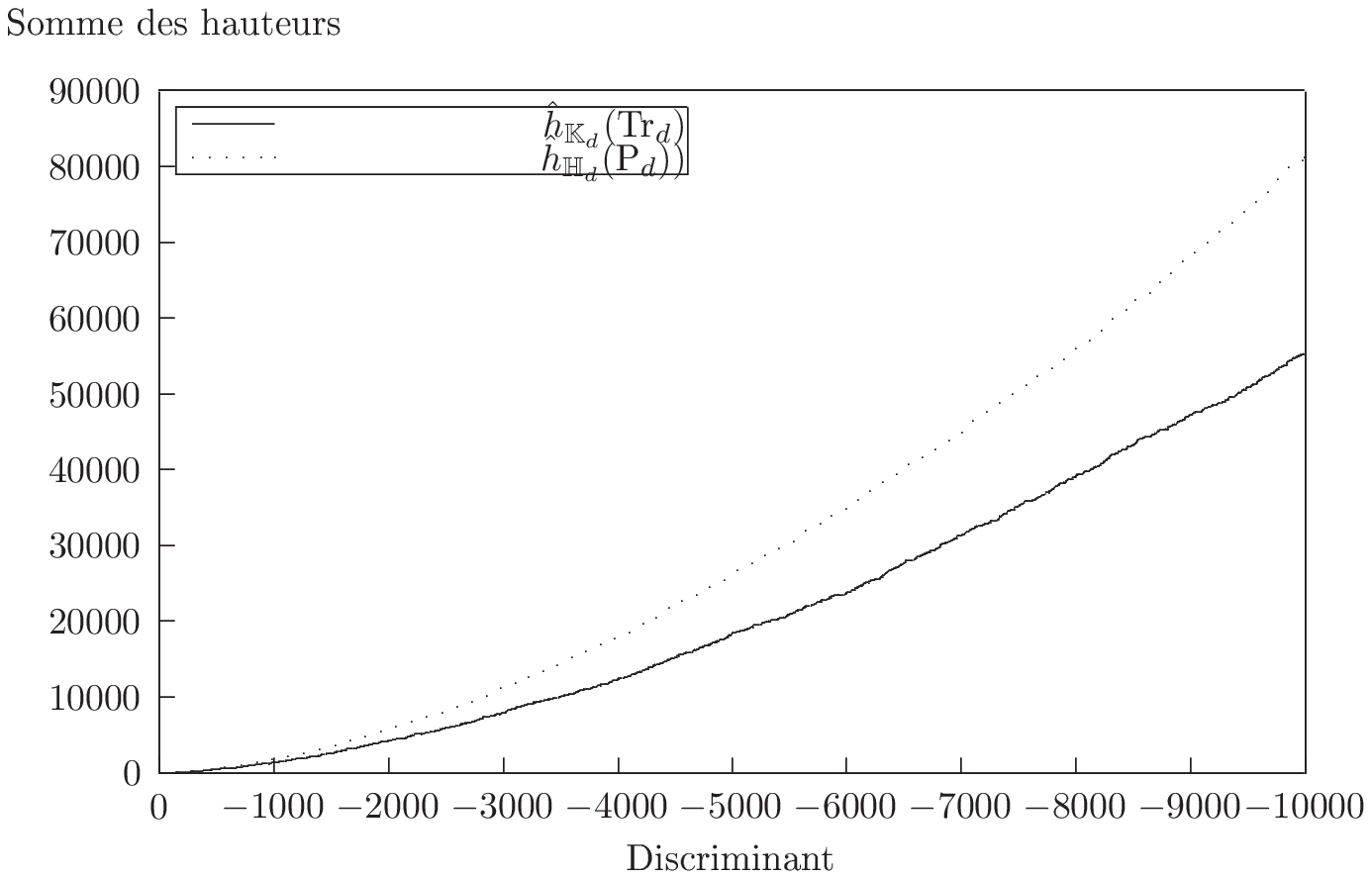}
\end{center}
\caption{Somme des hauteurs des traces et des points sur la courbe $37\text{B}$.}
\label{37B-points-traces}
\end{figure}

Pour \'etudier le rapport
$$\f{C_{\text{P}}}{C_{\text{Tr}}^{(0)}}= \f{\mathcal{Q}\mn{(N)}}{\mathcal{P}\mn{(1)}}L(E,1)^{-1}\prod_{p|N} \pa{1-\f{a_p}{p}} \pa{1-\f{a_{p^2}}{p^2}}^{-1},$$
on n\'eglige le r\^ole de $\f{\mathcal{Q}\mn{(N)}}{\mathcal{P}\mn{(1)}}$, qui est de toute fa\c{c}on born\'e. Ainsi, la taille de
$$\gamma_E =  L(E,1)^{-1} \prod_{p|N} \pa{1+\f{a_p}{p}}^{-1}$$
par rapport \`a $1$ refl\`ete essentiellement le signe du terme
 $$\sum_{\sigma\in G_d\backslash \{Id\}} <P,P^{\sigma}>_{\mathbb{H}_d}.$$
Plus $\gamma_E$ sera petit, plus ce produit scalaire sera grand, et, de mani\`ere imag\'ee, on pourrait dire que les points de Heegner sont essentiellement resser\'es autour d'une m\^eme direction ; alors que si $\gamma_E$ est grand devant $1$, cette somme est n\'egative et les points sont \'eclat\'es dans l'espace \`a $h_d$ dimensions. On s'attend donc, par exemple, \`a ce que la hauteur des traces (en moyenne) soit sup\'erieure \`a celle des points sur la courbe $58\text{B}$ ($\gamma_E=0.67$), ce qui est illustr\'e par la figure \ref{58B-points-traces}. Par contre, dans le cas de la courbe $37\text{B}$ ($\gamma_E=1.34$), les points sont plus gros (figure \ref{37B-points-traces}). 

\subsection{R\'esultats exp\'erimentaux}

Nous avons effectu\'e de nombreux calculs de points de Heegner et de leurs hauteurs sur diff\'erentes courbes, \`a l'aide des logiciels Magma et Pari. Magma a permis de calculer les points, ou les traces, eux-m\^emes (suivant la m\'ethode de calcul expos\'ee dans \cite{da}), alors que Pari s'est av\'er\'e plus rapide pour le calcul direct de la s\'erie $L$ intervenant dans la formule de Gross-Zagier. On s'est concentr\'e sur des courbes de petits conducteurs ($N<200$), car les algorithmes ont une complexit\'e en $\mathcal{O}(N^2)$. Nous pr\'esentons ici certains des r\'esultats obtenus, pour illustrer notre th\'eor\`eme.

\subsubsection{Comparaison entre les valeurs expérimentale et théorique de $C_{\text{P}}$}
Nous commen\c{c}ons par comparer les valeurs exp\'erimentales de $C_{\text{P}}$ \`a la valeur th\'eorique donn\'ee dans la section pr\'ec\'edente. Ainsi, pour chaque courbe de conducteur sans facteurs carr\'es plus petit que $100$, on a repr\'esent\'e le rapport entre la valeur exp\'erimentale
$$C_{\text{P}}^{\text{exp}}(Y):= \f{1}{Y^{3/2}\log Y}\sum_{\substack{d\in \mathcal{D} \\
\vert d\vert \leqslant Y}} \hat{h}_{\mathbb{H}_d}(\text{P}_d)$$
pour quelques valeurs de $Y$ et la valeur th\'eorique dans le tableau \ref{expe}.
\scriptsize
\begin{table}[htbp]
\begin{center}
$\begin{array}{ccccc}
\begin{array}{|c|c|c|c|}
\hline
\text{Courbe} & 6000 & 13000 & 20000 \\
\hline
11-1& 0.716 & 0.738 & 0.750 \\
14-1& 0.707 & 0.736 & 0.742 \\
15-1& 0.703 & 0.723 & 0.740 \\
17-1& 0.735 & 0.753 & 0.764 \\
19-1& 0.718 & 0.737 & 0.748 \\
21-1& 0.700 & 0.731 & 0.742 \\
26-1& 0.714 & 0.730 & 0.744 \\
26-2& 0.682 & 0.702 & 0.717 \\
30-1& 0.687 & 0.703 & 0.722 \\
33-1& 0.697 & 0.725 & 0.736 \\
34-1& 0.707 & 0.724 & 0.735 \\
35-1& 0.749 & 0.764 & 0.774 \\
37-1& 0.661 & 0.689 & 0.704 \\
37-2& 0.805 & 0.820 & 0.830 \\
38-1& 0.760 & 0.775 & 0.784 \\
38-2& 0.702 & 0.722 & 0.733 \\
39-1& 0.713 & 0.720 & 0.738 \\
\hline
\end{array}
& &
\begin{array}{|c|c|c|c|}
\hline
\text{Courbe} & 6000 & 13000 & 20000 \\
\hline
42-1& 0.683 & 0.722 & 0.735 \\
43-1& 0.678 & 0.699 & 0.714 \\
46-1& 0.699 & 0.720 & 0.729 \\
51-1& 0.711 & 0.719 & 0.737 \\
53-1& 0.690 & 0.716 & 0.725 \\
55-1& 0.818 & 0.821 & 0.825 \\
57-1& 0.680 & 0.695 & 0.705 \\
57-2& 0.774 & 0.780 & 0.785 \\
57-3& 0.727 & 0.738 & 0.745 \\
58-1& 0.663 & 0.687 & 0.703 \\
58-2& 0.808 & 0.818 & 0.828 \\
61-1& 0.735 & 0.752 & 0.770 \\
62-1& 0.832 & 0.835 & 0.847 \\
65-1& 0.752 & 0.756 & 0.766 \\
66-1& 0.752 & 0.773 & 0.785 \\
66-2& 0.693 & 0.719 & 0.734 \\
66-3& 0.671 & 0.699 & 0.714 \\
\hline
\end{array} 
& &
\begin{array}{|c|c|c|c|}
\hline
\text{Courbe} & 6000 & 13000 & 20000 \\
\hline
67-1& 0.725 & 0.753 & 0.763 \\
69-1& 0.795 & 0.797 & 0.807 \\
70-1& 0.762 & 0.774 & 0.785 \\
73-1& 0.803 & 0.822 & 0.834 \\
77-1& 0.709 & 0.738 & 0.752 \\
77-2& 0.767 & 0.791 & 0.802 \\
77-3& 0.770 & 0.795 & 0.806 \\
78-1& 0.722 & 0.722 & 0.742 \\
79-1& 0.794 & 0.823 & 0.827 \\
82-1& 0.731 & 0.754 & 0.763 \\
83-1& 0.786 & 0.799 & 0.811 \\
85-1& 0.808 & 0.820 & 0.826 \\
89-1& 0.823 & 0.838 & 0.844 \\
89-2& 0.800 & 0.817 & 0.824 \\
91-1& 0.738 & 0.748 & 0.754 \\
91-2& 0.766 & 0.772 & 0.778 \\
94-1& 0.996 & 0.980 & 0.974 \\
\hline
\end{array}\\
\end{array}$
\end{center}
\caption{Rapport entre valeur expérimentale $C_{\text{P}}^{\text{exp}}(Y)$ et valeur théorique $C_{\text{P}}$.}
\label{expe}
\end{table}
\normalsize
Ainsi, m\^eme pour des discriminants assez grands ($2\cdot 10^4$), la constante exp\'erimentale est souvent de l'ordre de $75\%$ de la constante th\'eorique. On a repr\'esent\'e plusieurs valeurs de $Y$ pour bien montrer que ce rapport augmente toutefois, mais tr\`es lentement.

\subsubsection{\'Etude plus fine des courbes $37\text{A}$ et $37\text{B}$}

Nous allons \'etudier plus en profondeur les courbes $37\text{A}$ et $37\text{B}$. Ces deux courbes sont int\'eressantes pour plusieurs raisons : elles ont m\^eme degr\'e et m\^eme conducteur, donc devraient avoir m\^eme $C_{\text{P}}$. La courbe $37\text{B}$ est de rang $0$ alors que la $37\text{A}$ est la courbe de rang $1$ de plus petit conducteur. 


\begin{figure}[htbp]
\begin{center}
\includegraphics[height=10cm]{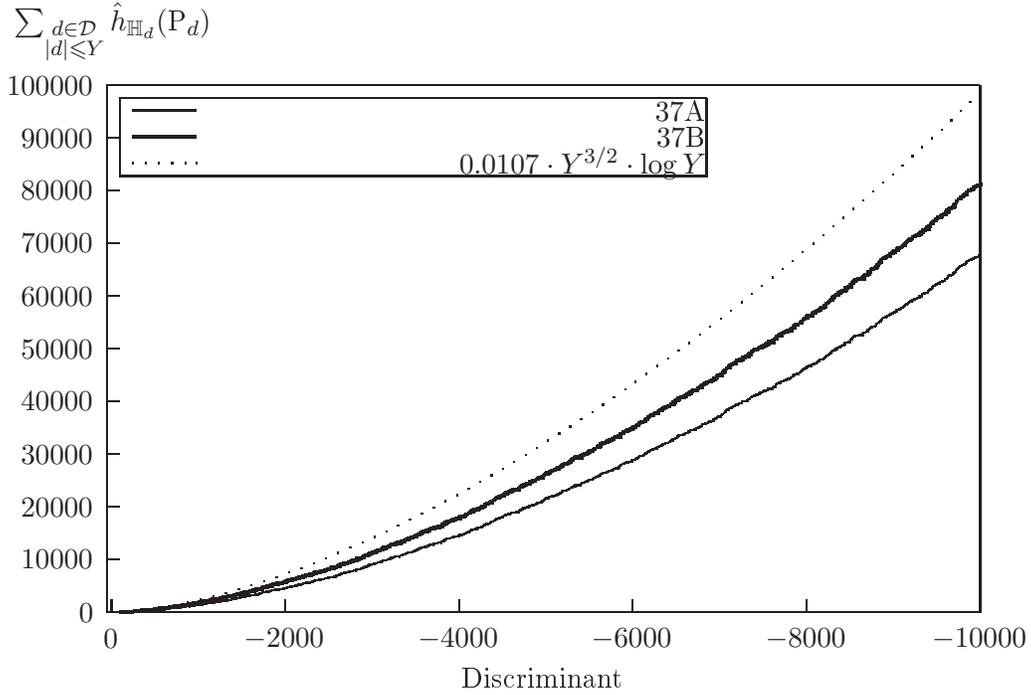}
\end{center}
\caption{Hauteur des points en moyenne sur les courbes $37\text{A}$ et $37\text{B}$.}
\label{37AB-points1}
\end{figure}



La figure \ref{37AB-points1} repr\'esente les sommes des hauteurs des points sur les courbes $37\text{A}$ et $37\text{B}$ compar\'ees \`a la valeur th\'eorique donn\'ee par le corollaire \ref{hpointstheo}. Contrairement \`a ce que l'on avait pour les hauteurs des traces, ici les courbes ne se supperposent pas du tout, ce qui \'etait pr\'evisible \'etant donn\'e le tableau ci-dessus.


On a en particulier l'impression que la courbe $37\text{A}$ est nettement en-dessous de la $37\text{B}$ sans para\^itre la rejoindre alors que le corollaire \ref{hpointstheo} affirme que les hauteurs des points sur ces courbes elliptiques devraient \^etre les m\^emes en moyenne. Cependant, une analyse plus fine de la diff\'erence entre ces deux courbes montre qu'elle semble \^etre en $Y^{3/2}$ et donc que les deux courbes semblent se rapprocher \`a une vitesse de $1/\log Y$ de la courbe théorique d'équation $Y\mapsto 0.0107Y^{\frac{3}{2}}\log{Y}$ ce qu'il est malheureusement difficile d'observer dans l'\'echelle de discriminants repr\'esent\'ee. Autrement dit, on devine numériquement sur les courbes $37$A et $37$B que
\begin{equation*}
\sum_{\substack{d\in\mathcal{D} \\
\vert d\vert\leqslant Y}}\hat{h}_{\mathbb{H}_d}(\text{P}_d)=C_{\text{P}}Y^{\frac{3}{2}}\log{Y}\left(1+\mathcal{O}_{N,E}\left(\frac{1}{\log{Y}}\right)\right).
\end{equation*}
Selon le corollaire \ref{hpointstheo} et la preuve du théorème \ref{traceana}, on a
\begin{eqnarray*}
\sum_{\substack{d\in\mathcal{D} \\
\vert d\vert\leqslant Y}}\hat{h}_{\mathbb{H}_d}(\text{P}_d) & = & C_{\text{P}}Y^{\frac{3}{2}}\log{Y}+C_{\text{P}}^\prime Y^{\frac{3}{2}}+\frac{1}{3\Omega_{E,N}}Y^{\frac{1}{2}}\mathsf{Error}+\mathcal{O}_{N,\varepsilon}\left(Y^{\frac{29}{20}+\varepsilon}\right), \\
& = & C_{\text{P}}Y^{\frac{3}{2}}\log{Y}\left(1+\frac{C_{\text{P}}^\prime}{C_{\text{P}}}\frac{1}{\log{Y}}+\frac{1}{3\Omega_{E,N}C_{\text{P}}}\frac{\mathsf{Error}}{Y\log{Y}}+\mathcal{O}_{N,\varepsilon}\left(Y^{-\frac{1}{20}+\varepsilon}\right)\right)
\end{eqnarray*}
où $\mathsf{Error}$ est défini en \eqref{err1} et pour tout $\varepsilon>0$. L'analyse numérique sugg\`ere donc que le terme $C_{\text{P}}Y^{3/2}$ dans le d\'eveloppement du corollaire \ref{hpointstheo} est non nul et m\^eme de l'ordre du terme principal pour des petits discriminants. Ceci suggère également que
\begin{equation}
\label{devine}
\mathsf{Error}=o_{\varepsilon}(NY).
\end{equation}
\label{fin}
Prouver cela nécessite de pouvoir estimer les moyennes mentionnées dans la remarque \ref{rq}. Donnons une autre justification numérique de nos intuitions. Posons
$$\delta(Y):=\f{1}{Y^{3/2}} \sum_{\substack{d\in\mathcal{D} \\
\vert d\vert\leqslant Y}} \pa{\hat{h}_{\mathbb{H}_d,37\text{B}}(\text{P}_d)-\hat{h}_{\mathbb{H}_d,37\text{A}}(\text{P}_d)}.$$
Le tableau \ref{delta} donne la valeur de $\delta(Y)$ pour plusieurs valeurs de $Y$.
\begin{table}[htbp]
\begin{center}
$\begin{array}{|c|c|c|c|c|c|c|}
\hline
Y  & 2\cdot 10^4    & 4\cdot 10^4    & 6\cdot 10^4   & 8\cdot 10^4  & 10\cdot 10^4 \\ 
\hline
\delta(Y)  & 0.01337  & 0.01329  & 0.01328  & 0.01326 & 0.01324 \\
\hline
\end{array}$
\end{center}
\caption{Valeurs numériques de $\delta(Y)$.}
\label{delta}
\end{table}
Ainsi, $\delta(Y)$ d\'ecro\^it tr\`es l\'eg\`erement avec $Y$ et semble se stabiliser. On devine alors que
\begin{equation*}
\delta(Y)=D_E+o_{N,E}(1)
\end{equation*}
pour une constante $D_E$. Or, le corollaire \ref{hpointstheo} affirme que
\begin{equation*}
\delta(Y)=\left(C_{\text{P},37\text{B}}-C_{\text{P},37\text{A}}\right)+\frac{1}{3Y}\left(\frac{\mathsf{Error}_{37\text{B}}}{\Omega_{37\text{B},37}}-\frac{\mathsf{Error}_{37\text{A}}}{\Omega_{37\text{A},37}}\right)+\mathcal{O}_{N,\varepsilon}\left(Y^{-\frac{1}{20}+\varepsilon}\right)
\end{equation*}
pour tout $\varepsilon>0$ ce qui confirme \eqref{devine}. En outre, il ne semble pas y avoir de compensation entre $C_{\text{P}}^\prime Y^{\frac{3}{2}}$ et $C_{\text{P}}^{\prime\prime}Y^{\frac{1}{2}}\mathsf{Err}_1$ car sinon $\delta(Y)$ tendrait plus vite vers $0$.

\strut\newline
\noindent{\textit{G. Ricotta}\newline}
Université de Montréal, Département de Mathématiques et de Statistique, CP 6128 succ Centre-Ville, Montréal QC   H3C 3J7, Canada; ricotta@dms.umontreal.ca\newline
\noindent{\textit{T. Vidick}\newline}
\'Ecole Normale Sup\'erieure, 45 rue d'Ulm, 75005 Paris, France; thomas.vidick@ens.fr

\end{document}